\documentclass[12pt]{article}
\usepackage{graphics,latexsym,amssymb,amsmath,amscd}
\RequirePackage{graphics}
\usepackage[all,cmtip]{xy}
\usepackage{xcolor}
\def\Im{\hbox{\rm Im}\,}
\def\[#1\]{\begin{eqnarray*}#1\end{eqnarray*}}

\def\tr{\hbox{\rm tr}\,}

\def\SU{\hbox{\bf SU}}
\def\SL{\hbox{\bf SL}}
\def\phi{\varphi}

\input xy
\xyoption{all}
\newtheorem{thm}{Theorem}[section]
\newtheorem{dfn}[thm]{Definition}

\newtheorem{cor}[thm]{Corollary}
\newtheorem{prop}[thm]{Proposition}
\newtheorem{lemma}[thm]{Lemma}

\newcommand{\Pf}{{\em Proof}. }
\newcommand{\EPf}{\hbox{}\hfill$\Box$\vspace{.5cm}}
\newcommand{\C}{{{\mathbb C}}}

\newcommand{\R}{{{\mathbb R}}}
\newcommand{\Z}{{{\mathbb Z}}}

\def\tr{\hbox{\rm tr}\,}

\def\SU{\hbox{\bf SU}}
\def\SL{\hbox{\bf SL}}
\def\PSL{\hbox{\bf PSL}}

\def\phi{\varphi}

\def\C{{\mathbb C}}
\def\R{{\mathbb R}}

\date{}
\title{A global invariant for path structures and second order differential equations}
\author{E. Falbel and J. M. Veloso}
 
\begin{document}

\maketitle
\newcommand{\D}{\mbox{$\cal D$}}
\newtheorem{df}{Definition}[section]
\newtheorem{te}{Theorem}[section]
\newtheorem{co}{Corollary}[section]
\newtheorem{po}{Proposition}[section]
\newtheorem{lem}{Lemma}[section]
\newcommand{\Ad}{\mbox{Ad}}
\newcommand{\ad}{\mbox{ad}}
\newcommand{\im}[1]{\mbox{\rm im\,$#1$}}
\newcommand{\bm}[1]{\mbox{\boldmath $#1$}}
\newcommand{\sime}{\mbox{sim}}
\begin{abstract}
We study a global invariant for path structures.  The invariant  is obtained as a secondary invariant from a Cartan connection on a canonical bundle associated to a path structure.   
It is  computed in examples which are defined in terms of reductions of the path structure.  In particular we give a formula for this global invariant for  second order differential equations defined on a torus $T^2$.

\end{abstract}

\section{Introduction}

Path structures on a 3-manifold are defined by a choice of contact structure and a decomposition of the contact plane bundle as a direct sum of two
line bundles.  This structure was througly studied in the 19th century (see in particular \cite{T}) as it appears in the description of 
second order differential equations and their equivalence under certain transformations (see Section \ref{section:definition} and references \cite{A,IL,BGH}).   

In Section \ref{section:definition} we collect definitions and examples. In particular we explain the relation with ordinary second order equations. In the following section we define the most important reductions of path structures.  The first one is obtained by fixing a global contact form and it is called strict path structure.  There exists a Cartan bundle $Y_1$ and  a connection adapted to that structure (see \ref{section:pseudo}) which was used in (\cite{FMMV}) to obtain a classification of compact 3-manifolds with non-compact automorphism group preserving the strict path structure.  We recall the construction in Proposition \ref{cartan-strict}. The second one, we call enriched path structure following \cite{MM} which were used  by Mion-Mouton  to classify certain classes of partial-hyperbolic diffeomorphisms of three manifolds.  It consists of path structures where we fix a line transverse to the contact distribution.  We define an adapted Cartan bundle $Y_2$ and a canonical connection adapted to this structure (see \ref{section:enriched} and Proposition \ref{Cartan-enriched}). 
 There exists a natural embedding $Y_1\to Y_2$ (Section \ref{subsectionY1toY2}, Proposition \ref{propY1toY2}).  

In Section \ref{section:Y} we recall the construction of the Cartan bundle $Y$ and the canonical adapted connection to a path structure on a 3-manifold (see Proposition \ref{Cartan-path}).  This construction is due to Cartan in \cite{Car}.  Although one can find modern treatments of this topic in several references (in particular \cite{IL, BGH}), we  include this section for the sake of completeness and because the conventions we use might differ from others.  We obtain a natural embedding $Y_2\to Y$ (see \ref{section:embedddingY2}, Proposition \ref{propY2toY}) and compute the curvature of the bundle $Y$ in terms of the curvature of $Y_2$ (see \ref{section:embedddingY2-curvature}).  The formulas are used in the computation of the global invariant in the next section.  We also recall the computations by Cartan of the invariants of a second order differential equation.

In the following section we define the global invariant when $Y_2$ admits a global section (see Definition \ref{definition:global}).  This construction is inspired by an analogous construction of a Chern-Simons invariant in the case of CR manifolds given in \cite{BE} (see also \cite{CL} for a relative version which does not depend  on the existence of a global section).  In \cite{FV} we defined a similar invariant for flag structures.  Those are manifolds equipped with a decomposition of a complex contact structure defined on the complexified tangent bundle of a 3-manifold.  In this paper we restrict the definition to path structures.  We obtain the expression of the invariant in terms of a reduction $Y_2$ or $Y_1$ of the Cartan bundle $Y$ of the path structure (see Proposition \ref{proposition-invariantY_2}). We also give a formula of the invariant in the case of a second order differential equation on the torus (Proposition \ref{proposition;global-invariant}).  It involves an integration of fifth order derivatives of  the function defining the ordinary equation in the form $y''=F(x,y,y')$.  We use coordinates in the projective cotangent bundle over a surface as explained in section \ref{section:globalDE}.  We characterize certain families of differential equations on the torus which have vanishing global invariant in Corollary \ref{corollary:Q1mu=0}. We then compute the invariant for a family of path structures on tight contact structures on the torus (see Proposition \ref{proposition:torustight}) and characterize those structures with vanishing global invariant, they turn to be flat. 
Finally we compute the global invariant for homogeneous path structures on $\SU(2)$ (see Proposition \ref{proposition:strictori}) and identify the  flat structure on the sphere where the global invariant is maximal. 

The authors thank Martin Mion-Mouton for useful discussions.

\section{Path structures in dimension 3}\label{section:definition}

Path geometries are very related to the theory of second order differential equations.  See  a modern  treatment in section 8.6 of \cite{IL} and in \cite{BGH} where the relation to second order differential equations is also explained.   
Le $M$ be a real three dimensional manifold and $TM$ be its tangent bundle.

\begin{dfn}A path structure on $M$ is a choice of two sub-bundles $T^1$ and $T^2$ in $TM$ such that 
$T^1 \cap T^2=\{ 0\}$ and such that $T^1 \oplus T^2$ is a contact distribution.
\end{dfn}

The condition that $T^1 \oplus T^2$ be a contact distribution means that, locally,  there exists a one form $\theta\in {T^*M}$ such that
$\ker \theta= T^1\oplus T^2$ and $d\theta \wedge \theta$ is never zero.

 One can choose a contact form $\theta$ up to a scalar function. One can interpret this as follows: one has a $\R^*$-bundle over the manifold given by the choice of $\theta$ at each point (one might keep only positive multiples for simplicity).  Over this line bundle one defines the tautological form $\omega_x=\pi^*(\theta_{\pi(x)})$.   
This bundle is trivial if and only if there exists a global contact form $\theta$.  {{If the contact distribution is oriented, then there exists a global contact form. Indeed, using a global metric on the distribution one can define locally a transversal vector to the distribution taking a Lie bracket of orthonormal vectors in the distribution. This defines a global 1-form.}}
 
Fix $\theta$ and local forms $Z^1$ and $Z^2$ defining the lines as above such that $d\theta =Z^1\wedge Z^2$.
There exists global forms $Z^1$ and $Z^2$ if and only if there exists global vector fields along the lines.  
Clearly, if the contact distribution is oriented, it suffices that there exists a global vector field along one of the foliations by lines.  

Local equivalence (also called point equivalence) between path structures happens when there exists a  local diffeomorphism which gives a correspondence between the lines defining each structure.

\subsection{The flat model space}
Flat path geometry   is the geometry of real flags in $\R^3$.  That is the geometry of the space of all couples $(p,l)$ where $p\in \R P^2$ and $l$ is a real projective line
containing $p$.  The space of flags is identified to the quotient
$$
\SL(3,\R)/B
$$
where $B$ is the Borel group of all real upper triangular matrices.  

The Lie algebra of $\SL(3,\R)$ decomposes into the following direct sum of vector subspaces:
$$
 {\mathfrak {sl}}(3,\R)= {\mathfrak {g}}^{-2}\oplus {\mathfrak {g}}^{-1}\oplus {\mathfrak {g}}^{0}\oplus{\mathfrak {g}}^{1}\oplus {\mathfrak {g}}^{2},
$$
where 
$${\mathfrak {g}}^{-2}= \left\{\left ( \begin{array}{ccc}

                        0      &    0    &   0   \\

                        0     &    0     & 0 \\

                        z       &  0    &     0
                \end{array} \right )\right\},\ \ \
 {\mathfrak {g}}^{-1}= \left\{\left ( \begin{array}{ccc}

                        0      &    0    &   0   \\

                        x     &    0     &   0 \\

                        0       &  y    &     0

                \end{array} \right )\right\},
$$ 
$$
{\mathfrak {g}}^{0}=\left\{\left ( \begin{array}{ccc}

                        u  +v    &    0    &  0   \\

                        0    &   -2v     &   0\\

                        0       &  0    &       -u+v

                \end{array} \right )\right\},
$$
$$
{\mathfrak {g}}^{1}=\left\{\left ( \begin{array}{ccc}

                        0      &    a    &   0   \\
                       0     &    0     &   b \\

                        0       &  0    &     0

                \end{array} \right )\right\},\ \ \
{\mathfrak {g}}^{2}=\left\{\left ( \begin{array}{ccc}

                        0      &    0    &   c   \\

                        0   &    0     &   0 \\

                        0       &  0    &     0

                \end{array} \right )\right\}.
$$

That is  the graded decomposition of  ${\mathfrak {sl}}(3,\R)$ where ${\mathfrak b}= {\mathfrak {g}}^{0}\oplus{\mathfrak {g}}^{1}\oplus {\mathfrak {g}}^{2}$ corresponds to upper triangular matrices with null trace.  The tangent space of $ \SL(3,\R)/B$ at $[B]$ is identified to
$$
{\mathfrak {sl}}(3,\R)/{\mathfrak b}= {\mathfrak {g}}^{-2}\oplus {\mathfrak {g}}^{-1}.
$$


\subsection{Examples}

\vspace{.5cm}
{\bf{Example I}} Consider the Heisenberg group 
$$
\mathbf{Heis}(3)=\{\ (z,t)\ \vert \ z\in \C,\, t \in \R\ \}
$$
with multiplication defined by $(z_1,t_1)\star (z_2,t_2)= (z_1+z_2, t_1+t_2+ 2\Im{z_1 \overline{z_2}})$.  The contact form 
$$
\theta= dt-xdy-ydx
$$
is invariant under left multiplications (also called Heisenberg translations).  If $\Lambda\subset \mathbf{Heis}(3)$ is a lattice then the quotient
$\Lambda\setminus \mathbf{Heis}(3)$ is a circle bundle over the torus with a globaly defined 
contact form.  

A lattice $\Lambda$ determines a lattice  $\Gamma\subset \C$ corresponding to projection in the exact sequence
$$
0\rightarrow \R\rightarrow \mathbf{Heis}(3)\rightarrow \C\rightarrow 0.
$$
There are many global vector fields in the distribution defined by $\theta$ invariant under $\Lambda$, it suffices to lift an invariant vector field on $\C$ under $\Gamma$.
All circle bundles obtained in this way are not trivial and the fibers are 
transverse to the distribution.

{\bf{Example II}}. Here we consider the torus $T^3$ with coordinates $(x,y,t)$ ($\mod 1$) and the global contact form
$$
\theta_n= \cos (2\pi n t) dx- \sin(2\pi n t) dy.
$$
There are two canonical global vector fields on the distribution given by
$\frac{\partial}{\partial t}$ and $\sin (2\pi n t) \frac{\partial}{\partial x}+\cos(2\pi n t)\frac{\partial}{\partial y}$.
In this example, the fiber given by the coordinate $t$ has tangent space contained in the distribution.

{\bf{Example III}}.  An homogeneous example is the Lie group $\SU(2)$ with  left invariant vector fields $X$ and $Y$ with $Z=[X,Y]$ and cyclic commutation relations.
The vector fields $X$ and $Y$ define a path structure on $\SU(2)$.

{\bf{Example IV}}.  Another  homogeneous example is the Lie group $\SL(2,\R)$ with  left invariant vector fields $X$ and $Y$ with $Z=[X,Y]$ with $[Z,X]=X$ and $[Z,Y]=-Y$ given by generators 
$$
X= \left ( \begin{array}{cc}

                        0      &    1    \\

                        0 &     0
                                             
                \end{array} \right ), \ 
                Y= \left ( \begin{array}{cc}

                        0       &    0    \\

                        1 &     0
                                             
                \end{array} \right ), \ 
                Z= \left ( \begin{array}{cc}

                        1       &    0    \\

                        0 &     -1
                                             
                \end{array} \right ).
$$
The path structure defined by $X$ and $Y$ induces a path structure on the quotient $\Gamma\setminus\SL(2,\R)$ by a discrete torsion free subgroup $\Gamma\subset \SL(2,\R)$.  This structure is invariant under the flow defined by right multiplication by $e^{tZ}$.

{\bf{Example V}}.  Let $\Sigma$ be a surface equipped with a Riemannian metric.   The geodesic flow on the unit tangent bundle $T^1\Sigma$ defines a distribution which, together with the distribution defined by the vertical fibers of the projection of the unit tangent bundle on $\Sigma$, defines a path structure which is not invariant under the geodesic flow.  For $\Sigma =H^2_\R$, the hyperbolic space, we obtain  $T^1\Sigma= \PSL(2,\R)$ with distributions defined by the left invariant distributions
$X-Y$ and $Z$ (using the same generators of the Lie algebra as in the previous example).

{\bf{Example VI}}
Let $M$ be a three manifold equipped with a  path structure $D=T^1\oplus T^2\subset TM$.  Suppose  $D$ is orientable and choose a section $u$ of $T^1$.  Each section $v$ of $T^2$ such that $(u,v)$ is positive gives rise to a CR structure.  Indeed we define $Ju=v$ and $Jv=-u$.   The choice of $v$ corresponds to a section of an $\R^*_+$-bundle over $M$.
Reciprocally given a CR structure on $M$, defined by $J:D\rightarrow D$, one can associate path structures corresponding to a choice $T^1\subset D$ and defining then $T^2=J(T^1)$.

\subsection{Path structures and second order differential equations}

This is studied since a long time (see \cite{T}, \cite{IL} and \cite{BGH}).  It turns out that path structures can be obtained putting together second order differential equations in one variable.  Indeed, a second order differential equation in one variable is described locally as
$$
\frac{d^2 y}{dx^2}=F(x,y,\frac{dy}{dx}).
$$
This defines a path structure on a neighborhood of a point in $\R^3$  with coordinates $(x,y,p)$:
$$
L_1= \ker\{d p-Fdx\}\cap \ker\{dy-pdx\},\ \ \ L_2= \ker dx\cap \ker dy.
$$
The contact structure is defined by the form
$$
\theta= d y-pdx.
$$
Defining the forms $Z^1=dx$ and $Z^2=d p-Fdx$, one has that $d\theta= Z^1\wedge Z^2$.

One can show easily that every path structure is, in fact, locally equivalent to a second order equation.  That is, there exists local coordinates such that $L_1$ and $L_2$ are defined via a second order ODE as above.

\subsection{Reductions of path structures}\label{section:reductions}

 We will describe two reductions of path geometry corresponding to  subgroups $G_1\subset G_2\subset SL(3,\R)$ where
 $$
 G_1=\left\{\left ( \begin{array}{ccc}

                        a     &    0    &   0   \\
                       \star     &    \frac{1}{a^2}    &   0\\

                        \star       &  \star   &     a

                \end{array} \right )\right\}\ \ \
$$ and 
$$
 G_2=\left\{\left ( \begin{array}{ccc}

                        a     &    0    &   0   \\
                       \star     &    \frac{1}{ab}    &   0\\

                        \star       &  \star   &     b

                \end{array} \right )\right\}.
$$
The models are $G_1/\R^*$ and $G_2/{\R^*}^2$ and correspond to the Heisenberg group where in the first model we fix a contact form and, in the second, a transverse line to the contact distribution.

Other reductions of the $G_2$-structure might occur, namely by choosing other 
embeddings of $\R^*$ into $G_2$.  They appear naturally when certain components of the curvature 
of  the Cartan connections  on $Y_2$ or $Y$ are non-vanishing.

We will construct coframe bundles $Y_1, Y_2$ and a principal bundle $Y$ over $M$ with structure groups $\R^*, {\R^*}^2$ and the Borel group $B$ together with Cartan connections and canonical embeddings
$$
Y_1\rightarrow Y_2\rightarrow Y.
$$
They correspond to a strict path structure, an enriched path structure (see next sections for definitions) and finally, a path structures on the manifold $M$.

 \subsection{Path structures with a fixed contact form: strict path structures.}\label{section:pseudo}
 
 In this section we fix a contact form and recall the reduction of the structure group of a path geometry obtained in \cite{FV} where we called the path structure with a fixed contact form a pseudo flag structure.  This structure is called strict path structure in \cite{FMMV}.

 $G_1$ denotes from now on the subgroup of $\SL(3,\R)$ defined by
 $$
 G_1=
 \left\{\left ( \begin{array}{ccc}

                        a     &    0    &   0   \\
                       x    &    \frac{1}{a^2}    &   0\\

                        z       &  y &    a

                \end{array} \right )\ \vert\ {a\in\R^*,(x,y,z)\in\R^3}\right\}
$$ 
and $P_1\subset G_1$ the subgroup defined by
$$
 P_1=\left\{\left ( \begin{array}{ccc}

                        a     &    0    &   0   \\
                       0    &    \frac{1}{a^2}    &   0\\

                        0      &     0 &     a

                \end{array} \right )\right\}.\ \ \
$$
 Writing the Maurer-Cartan form of $G_1$ as the matrix
$$
\left ( \begin{array}{ccc}

                        w     &    0    &   0   \\
                       \theta^1    &    -2w   &   0\\

                       \theta      &  \theta^2   &     w

                \end{array} \right )
$$
one obtains the Maurer-Cartan equations:
$$d \theta+\theta^2\wedge \theta^1=0$$
$$d\theta^1 -3w\wedge \theta^1=0$$
$$d\theta^2 +3w\wedge \theta^2=0$$
$$d w=0.$$

{
 $G_1$ is the automorphism group of the canonical left-invariant strict path structure of $\mathbf{Heis}(3)$,
and that its action induces an identification of $\mathbf{Heis}(3)$ with the homogeneous space $X=G_1/P_1$.
}

Let $M$ be a three-manifold equipped with a strict path structure $(E^1,E^2,\theta)$ defined by two one dimensional bundles $E^1$ and $E^2$ and contact form $\theta$.
We let $R$ be the associated Reeb vector field (satisfying $\iota_Rd\theta=0$ and $\theta(R)=1$).
Now let  $X_1\in E^1$, $X_2\in E^2$ be such that $d\theta(X_1,X_2)=1$.  
The dual coframe of $(X_1,X_2,R)$ is $(\theta^1,\theta^2,\theta)$,
for two 1-forms $\theta1$ and $\theta_2$ verifying $d\theta=\theta^1\wedge \theta^2$.   

At any point $x\in M$,
one can look at the coframes of the form
$$
\omega^1={a^3}\theta^1(x), \  \omega^2=\frac{1}{a^3}\theta^2(x), \ \omega=\theta(x)
$$
for $a\in\R^*$.
\begin{dfn}
We denote by $p_1:Y_1 \rightarrow M$ the  {$\R^*$}-coframe bundle over $M$ 
  given by the set of coframes  $(\omega, \omega^{1}, \omega^{2})$ of the above form.
\end{dfn}

We will denote the tautological forms defined by $\omega^1,\omega^2,\omega$ using the same letters. 
That is, we write $\omega^i$ at the coframe $(\omega^1,\omega^2,\omega)$ to be $p_1^*(\omega^i)$.

\begin{prop}\label{cartan-strict}
There exists a unique Cartan connection on $Y_1$ 
$$
\pi_1=\left ( \begin{array}{ccc}

                        w     &    0    &   0   \\
                       \omega^1    &    -2w   &   0\\

                       \omega      &  \omega^2   &     w

                \end{array} \right )
$$ such that its curvature form is of the form
$$
\Pi_1=d\pi_1+\pi_1\wedge \pi_1=\left ( \begin{array}{ccc}

                        dw     &    0    &   0   \\
                       \omega\wedge \tau^1    &    -2dw   &   0\\

                       0      &  -\omega\wedge \tau^2   &    dw

                \end{array} \right )
$$
with $\tau^1\wedge \omega^2=\tau^2\wedge \omega^1=0$.
\end{prop}

Observe that the condition $\tau^1\wedge \omega^2=\tau^2\wedge \omega^1=0$ implies that we may write
$\tau^1=\tau^1_2\omega^2$ and $\tau^2=\tau^2_1\omega^1$.  The structure equations are
$$d \omega+\omega^2\wedge \omega^1=0,$$
$$d\omega^1 -3w\wedge \omega^1=\omega\wedge \tau^1,$$
$$d\omega^2 +3w\wedge \omega^2=-\omega\wedge \tau^2.$$

The proof of the proposition is given in \cite{FMMV} and \cite{FV}.

Bianchi identities are obtained differentiating the structure equations. They are described in the following equations:

 \begin{equation}\label{dtheta11}
dw=C\omega\wedge\omega^1+ D\omega\wedge\omega^2+S\omega^1\wedge\omega^2,
\end{equation}
\begin{equation}\label{dtau1}
d\tau^1_2-6\tau^1_2w+3D\omega^1=\tau^1_{20}\omega+\tau^1_{22}\omega^2
\end{equation}
\begin{equation}\label{dtau2}
d\tau^2_1+6\tau^2_1w+3C\omega^2=\tau^2_{10}\omega+\tau^2_{11}\omega^1
\end{equation}

\subsection{Path structures with a fixed transverse line: enriched path structures.}\label{section:enriched}

In this section we introduce a coframe bundle and a Cartan connection associated 
to a path structure with a fixed transverse line to the the contact distribution.

The model space is the homogeneous space which is the quotient of the group of lower triangular matrices in $SL(3,\R)$ by the subgroup of diagonal matrices.  The Maurer-Cartan form
is the Lie algebra valued form which can be represented by

$$\pi=\left ( \begin{array}{ccc}
                      \phi+w     &   0   &   0  \\

                     \omega^1          &   -2w   &    0 \\

                      \omega         &  \omega^{2}   &     -\phi+w
                \end{array} \right )
$$
The Maurer-Cartan equations $d\pi+\pi\wedge\pi=0$ are given by
$$
d\omega=2\phi \wedge \omega + \omega^1\wedge\omega^2
$$
$$
d\omega^1=\phi \wedge  \omega^1+3w\wedge \omega^1
$$
$$
d\omega^2=\phi\wedge \omega^2-3w\wedge \omega^2.
$$

Let $M$ be a three manifold equipped with a  path structure $D=T^1\oplus T^2\subset TM$. 
We fix a transverse line $L$ so that $TM= T^1\oplus T^2\oplus L$.

We suppose $X_1\in T^1$, $X_2\in T^2$ and $X\in L$ is a frame.  The dual coframe is $\theta^1$, $\theta^2$
and $\theta$.  Observe that $\theta$ is simply a form with $\ker \theta= D$. One can define a  coframe bundle defined  by all coframes:
$$
\omega^1=a^1\theta^1, \  \omega^2=a^2\theta^2, \ \omega=\lambda\theta.
$$
where we will suppose, for simplicity, that $a^1,a^2,\lambda>0$.

A reduction of this coframe bundle is obtained by imposing that each coframe verifies
$$
d\omega_{\vert D} = (\omega^1\wedge \omega^2) _{\vert D}
$$
for an extension of the 1-form such that $\ker \omega =D$.
This relation does not depend on the particular extension of a form $\omega$ defined at a 
point because $d \omega_{\vert D}(X,Y)=-\omega([X,Y])$ for any vector fields $X$ and $Y$ which are sections of the distribution $D$.

\begin{dfn}
We denote by $p_2: Y_2\rightarrow M$ the ${\R^*}^2$-coframe bundle over $M$ 
  given by the set of 1-forms  $(\omega, \omega^{1},
  \omega^{2})$ defined above.
The structure group ${\R^*}^2$ acts as follows
$$
(\omega', \omega'^{1},
  \omega'^{2})=
(\omega, \omega^{1},
  \omega^{2})
 \left ( \begin{array}{cccc}

                        \lambda       &    0               	&       0 \\

                        0  	&      a^1 		& 0             \\
                        0   	&  0 			&  a^2   \\

                \end{array} \right )
$$
where   $\lambda, a^1,a^2\in \R _+^*$ with $a^1 a^2=\lambda$.
\end{dfn}

In order to define a Cartan connection on $Y_2$ we start taking the 
tautological forms corresponding to the forms $\omega, \omega^1, \omega^2$, which we will denote by the same letters by abuse of notation.

Using a coframe section $(\theta, \theta^1,\theta^2)$ on $M$ one can express the tautological
forms as
$$
\omega =\lambda p_2^*(\theta), \ \omega^1=a^1p_2^*(\theta^1), \  \omega^2=a^2p_2^*(\theta^2),
$$
with $a^1a^2=\lambda$.

We need to define two forms $\phi$ and $w$ corresponding to the vertical directions 

Observe first that we have
$$
d\omega= \frac{d\lambda}{\lambda} \wedge \omega+\omega^1\wedge \omega^2 \ \ \mbox{mod}(\omega)
$$
and therefore one may write
\begin{equation}\label{do}
d\omega= 2\phi \wedge \omega+\omega^1\wedge \omega^2
\end{equation} 
where $\phi$ restricted to the vertical fiber is $\frac{d\lambda}{2\lambda}$.  The form $\phi$
is not yet fixed and any other form $\phi'$ satisfying the equation satisfies
$$
\phi-\phi'=s\omega
$$
where $s$ is a funtion on $Y_2$.

Differentiating the forms $\omega^1$ and $\omega^2$ we obtain new forms which correspond to the coordinates $a^1, a^2$ :

$d\omega^1= \frac{d a^1}{a^1}\wedge \omega^1 + a^1 d\theta^1$ and $d\omega^2= \frac{d a^2}{a^2}\wedge \omega^2+ a^2 d\theta^2$.  

Observing  that 
$$
\frac{d\lambda}{\lambda}= \frac{da^1}{a^1}+\frac{da^2}{a^2}
$$
we can write
$$
d\omega^1= \frac{d\lambda}{2\lambda} \wedge \omega^1+ \frac{1}{2}\left(\frac{da^1}{a^1}-\frac{da^2}{a^2}\right)\wedge \omega^1 
+a^1d \theta^1$$
$$
d\omega^2 = \frac{d\lambda}{2\lambda}\wedge  \omega^2- \frac{1}{2}\left(\frac{da^1}{a^1}-
\frac{da^2}{a^2}\right)\wedge \omega^2+a^2d \theta^2
$$

Now we can make the first right hand term of each equation to be $ \phi \wedge \omega^1$ and $ \phi \wedge \omega^2$ respectively by adding terms in $\omega, \omega^1, \omega^2$ to $\frac{d\lambda}{2\lambda}$.   The terms in $\omega^1\wedge \omega^2$  not appearing in these first terms can be absorbed in the second term in each equation.  It remains a last term in each equation that we denote by $\omega \wedge \tau^1$ and $-\omega\wedge \tau^2$ respectively.  We proved the following:

\begin{lem}
There exists  forms $w, \tau^1,\tau^2$ defined on $Y_2$
such that
\begin{equation}\label{domegaitau}
d\omega^1=\phi \wedge \omega^1+
3 w\wedge \omega^1 +\omega \wedge \tau^1\  \mbox{and}\ \ 
d\omega^2=\phi \wedge \omega^2-
3w\wedge \omega^2-\omega\wedge \tau^2.
\end{equation}
The forms $\tau^1$ and $\tau^2$ are horizontal, that is, they vanish on vectors tangent to the fibers of $Y_2\rightarrow M$.  Moreover,   writing $\omega^1= a^1\theta^1$, $\omega^2= a^2\theta^2$, 
$\omega=\lambda \theta$ for a choice of sections on $M$, one has   $\phi = \frac{d\lambda}{2\lambda}$ and $6w=\frac{da^1}{a^1}-\frac{da^2}{a^2}$ modulo the tautological forms of the fiber bundle $Y_2$. 

\end{lem}

Let  $\phi', w', {{\tau'}^1}$ and ${\tau'^2}$ be other forms satisfying
equations above. Taking the difference
we obtain
$$
0=(\phi-\phi')\wedge \omega^1+3(w-w')\wedge \omega^1+\omega\wedge (\tau^1-\tau'^1)
$$
and
$$
0=(\phi-\phi')\wedge \omega^2 -3(w-w')\wedge \omega^2 -\omega\wedge (\tau^2-\tau'^2)
$$
Therefore, as $\phi-\phi'=s\omega$, we can write
$$
0=-3\omega^1\wedge(w-w') +\omega\wedge (s \omega^1+\tau^1-\tau'^1)
$$
and
$$
0=3\omega^2\wedge(w-w') -\omega\wedge (-s \omega^2+\tau^2-\tau'^2).
$$
By Cartan's lemma we obtain
$$
w-w' = a\omega,
$$
$$
\tau^1-\tau'^1 = -3a\omega^1 -{s}\omega^1+b^1\omega,
$$
$$
\tau^2-\tau'^2 = -3a\omega^2 +{s}\omega^2+b^2\omega.
$$  

Now, we can impose that $\tau^1=\tau^1_1\omega^1+\tau^1_2\omega^2$ and $\tau^2=\tau^2_1\omega^1+\tau^2_2\omega^2$ by choosing convenient $b^1$ and $b^2$ (or by simply considering, from the beginning, $\tau^1$ and $\tau^2$ with no terms in $\omega$).  Moreover, one can choose  unique $a$ and $s$  so that $\tau^1_1=0$ and  $\tau^2_2=0$.  We conclude that

\begin{lem}
There exists  unique forms $\phi, w, \tau^1,\tau^2$ defined on $Y_2$
such that
$$
d\omega=2\phi \wedge \omega + \omega^1\wedge\omega^2
$$
$$
d\omega^1=\phi \wedge  \omega^1+3w\wedge \omega^1 +\omega \wedge \tau^1
$$
$$
d\omega^2=\phi\wedge \omega^2-3w\wedge \omega^2 -\omega\wedge \tau^2
$$
with  $\tau^1\wedge \omega^2=\tau^2\wedge \omega^1=0$.
\end{lem}

Bianchi identities are obtained differentiating the above equations:

\begin{lem}
There exists a 1-form $\psi$ such that
\begin{equation}\label{dphi}
d\phi = 
\omega \wedge \psi
\end{equation}
The form $\psi$ may be chosen satisfying $\psi= A\omega^1+B\omega^2$ and $d\psi =-2\phi\wedge \psi +\omega\wedge \alpha$ where $A, B$ are functions on $Y_2$ and $\alpha$ is a 1-form on $Y_2$.
\end{lem}
\Pf
Differentiating equation $d\omega=2\phi \wedge \omega + \omega^1\wedge\omega^2$ one obtains, using equations \ref{domegaitau},
that $d\phi\wedge \omega=0$, that is,
\begin{equation}
d\phi = 
\omega \wedge \psi
\end{equation}
for a 1-form $\psi$ defined on $Y_2$.

Differentiating  $d\phi = \omega \wedge \psi$ one has
$$
0=d\omega\wedge \psi -\omega\wedge d\psi=(2\phi \wedge \omega + \omega^1\wedge\omega^2)\wedge \psi -\omega\wedge d\psi=\omega^1\wedge\omega^2\wedge \psi-\omega\wedge(d\psi+2\phi\wedge \psi).
$$
Using Cartan's lemma,  $\psi=A\omega^1+B\omega^2$ modulo $\omega$, and we certainly can choose $\psi$ satisfying 
$d\phi = 
\omega \wedge \psi$ with $\psi=A\omega^1+B\omega^2$.  We conclude that
$$
d\psi+2\phi\wedge \psi=\omega\wedge \alpha.
$$
\EPf

Equation $dd\omega^1=0$ 
gives after simplifications

\begin{align}\label{ddomega1}
0=d(\phi+3w)\wedge \omega^1+ \omega\wedge\omega^2 (d\tau^1_2+2\tau^1_2(\phi-3w)).
\end{align}

Analogously, $dd\omega^2=0$  
simplifies to

\begin{align}\label{ddomega2}
0=d(\phi-3w)\wedge \omega^2- \omega\wedge\omega^1 (d\tau^2_1+2\tau^2_1(\phi+3w)).
\end{align}

Using the previous lemma 
we may write
$$
dw=C\omega\wedge\omega^1+ D\omega\wedge\omega^2+S\omega^1\wedge\omega^2,
$$
where $C,D$ and S are functions on $Y_2$.

We can represent the equations above 
  as a matrix equation whose entries are differential forms. The forms are disposed in the Lie algebra ${\mathfrak {b}}\subset {\mathfrak {sl}}(3,\R)$ (the Lie algebra of lower triangular matrices) and we obtain the following Proposition. 

\begin{prop}\label{Cartan-enriched}

Let $Y_2$ be the adapted principal bundle constructed above associated to an enriched path structure on a manifold $M$.  Then there exists a unique Cartan's connection with values in ${\mathfrak {b}}$:  
$$\pi_2=\left ( \begin{array}{ccc}
                        \phi+w     &   0   &   0  \\

                     \omega^1          &   -2w   &    0 \\

                      \omega         &  \omega^{2}   &      -\phi+w
                \end{array} \right )
$$
with curvature:
\begin{equation}\label{dpi}
\Pi_2=d\pi_2+\pi_2\wedge\pi_2=\left ( \begin{array}{ccc}
                       \omega\wedge \psi  +W&   0      &  0   \\

                         \tau^1_2\omega\wedge \omega^2  &   -2W    &  0  \\

                           0  &  - \tau^2_1\omega\wedge \omega^1  &     -\omega\wedge \psi +W
               \end{array} \right )
\end{equation}
where $
W=C\omega\wedge \omega^1+D\omega\wedge \omega^2+ S\omega^1\wedge \omega^2
$ and $\psi= A \omega^1+B \omega^2$.
\end{prop}

\subsubsection{More Bianchi identities}

\begin{itemize}

\item
Substituting the expressions above in equations \ref{ddomega1} and \ref{ddomega2} we obtain
\begin{equation}\label{dtau12}
d\tau^1_2+2\tau^1_2(\phi-3w)+(B+3D)\omega^1=\tau^1_{20}\omega+\tau^1_{22}\omega^2.
\end{equation}

\item 
Analogously we obtain
\begin{equation}\label{dtau21}
d\tau^2_1+2\tau^2_1(\phi+3w)-(A-3C)\omega^2=\tau^2_{10}\omega+\tau^2_{11}\omega^1.
\end{equation}

From the last two equations we obtain the following

\begin{prop} 
\label{prop:nultorsion}
If the adapted connection of $Y_2$ has nul torsion and 
$$
dw=S\omega^1\wedge \omega^2,
$$
then $d\phi=0$.
\end{prop}

\item Analogously, $dd\phi=0$  
simplifies to
$$
0=\omega\wedge\omega^1(dA+3A(\phi +w))+\omega\wedge\omega^2(dB+3B(\phi -w))
$$
and we obtain 
\begin{align}\label{dA}
dA+3A(\phi +w)=A_0\omega+A_1\omega^1+A_2\omega^2,
\end{align}
\begin{align}\label{dB}
dB+3B(\phi -w)=B_0\omega+B_1\omega^1+B_2\omega^2,
\end{align}
with $A_2=B_1$.

\item Also, $ddw=0$  
simplifies to
$$
0=\omega\wedge\omega^1(dC+3C(\phi +w))+\omega\wedge\omega^2(dD+3D(\phi -w))+\omega^1\wedge\omega^2(dS+2S\phi )
$$
and we obtain 
\begin{align}\label{dC}
dC+3C(\phi +w)=C_0\omega+C_1\omega^1+C_2\omega^2,
\end{align}
\begin{align}\label{dD}
dD+3D(\phi -w)=D_0\omega+D_1\omega^1+D_2\omega^2,
\end{align}
\begin{align}\label{dS}
dS+2S\phi=S_0\omega+S_1\omega^1+S_2\omega^2,
\end{align}
with $C_2-D_1+S_0=0$.

\end{itemize}

\begin{lem}

If $\tau^1=\tau^2=C=D=0$ 
$$
d\phi=0.
$$

\end{lem}

\Pf From the last formulae we obtain that $\psi$ is a multiple of $\omega$ and the result follows.
\EPf


\subsubsection{The embedding  $\iota_1: Y_1\to Y_2$}\label{subsectionY1toY2}

Given a path structure with a fixed contact form $\omega$ we obtained first a coframe bundle $Y_1$ and one can also obtain a canonical transverse direction by considering 
the Reeb vector field associated to $\omega$.  One obtains then a coframe bundle $Y_2$ of last section.

Given a coframe $(\omega, \omega^1,\omega^2)\in Y_1$ one can view the same coframe as a coframe of $Y_2$.  This gives an embedding
$$
\iota_1: Y_1\rightarrow Y_2.
$$

By abuse of language we may write the connection forms of $Y_1$ and $Y_2$ using the same letters and then obtain: 

\begin{prop}\label{propY1toY2}
There exists a unique embedding  $\iota_1: Y_1\rightarrow Y_2$ satisfying
$\iota_1^*(\omega)=\omega$, $\iota_1^*(\omega^1)=\omega^1$ and $\iota_1^*(\omega^2)=\omega^2$.  Moreover, for this embedding, $\iota_1^*(\phi)=0$ and  $\iota_1^*(w)=w$.
\end{prop}
\Pf If unicity is not satisfied one can obtain the same forms pulling back a different coframe.  But from the transformations of the coframe, 
\begin{align}
\tilde\omega  &= \frac{a}{b}\, \omega\notag\\
\tilde{\omega}^1 &=a^2b\, \omega^1\notag\\
\tilde\omega^2 &=\frac{1}{ab^2}\,\omega^2\notag.
\end{align} 
We obtain then that $a=b=1$ and the embedding is uniquely determined by the conditions.

Comparing the structure equations of both structures we further get the equations $\iota_1^*(\phi)=0$ and  $\iota_1^*(w)=w$.
\EPf

\section{The Cartan connection of a path structure}\label{section:Y}

We review in this section the construction of a Cartan connection.   The construction is due to E. Cartan in \cite{Car} and one can read a modern description 
of it in \cite{IL}.  We include this section in order to fix our conventions and describe the embedding of $Y_2$ into the corresponding fiber bundle associated to a path geometry (see \ref{section:embedddingY2} and \ref{section:embedddingY2-curvature}) which will be used to define the global invariant.

 The Maurer-Cartan form on $SL(3,\R)$ is given by a form  with values  in the Lie algebra ${\mathfrak {sl}}(3,\R)$ :

$$\pi=\left ( \begin{array}{ccc}
                        \phi+w     &   \phi^2   &   \psi  \\

                     \omega^1          &   -2w   &    \phi^1 \\

                      \omega         &  \omega^{2}   &      -\phi+w
                \end{array} \right )
$$
satisfying the equation $d\pi+\pi\wedge \pi =0$.  That is 
$$
d\omega = \omega^1\wedge \omega^2 +2\phi\wedge \omega
$$
$$
d\omega^1=\phi \wedge  \omega^1+3w\wedge \omega^1 +\omega \wedge \phi^1
$$
$$
d\omega^2=\phi\wedge \omega^2-3w\wedge \omega^2 -\omega\wedge \phi^2
$$
$$
dw= -\frac{1}{2}\phi^2\wedge \omega^1+\frac{1}{2}\phi^1\wedge \omega^2
$$
$$
d\phi = \omega \wedge \psi -\frac{1}{2}\phi^2\wedge \omega^1-\frac{1}{2}\phi^1\wedge \omega^2
$$
$$
d\phi^1 = \psi\wedge \omega^1-\phi\wedge \phi^1+3w\wedge \phi^1
$$
$$
d\phi^2 = -\psi\wedge \omega^2-\phi\wedge \phi^2-3w\wedge \phi^2
$$
$$
d\psi= \phi^1\wedge \phi^2+2\psi\wedge \phi.
$$

\subsection{The $\R^*$-bundle of contact forms and an adapted coframe bundle}\label{section:coframes}

We recall the construction of  the $\R^*$-bundle of contact forms.  
Define ${E}$ to be the $\R^*$-bundle of all forms $\theta$ on $TM$ such that $\ker \theta=T^1\oplus T^2$.   Remark that this bundle is trivial if and only if there   exists a globally defined 
non-vanishing form $\theta$. Define the set of forms $\theta^1$ and $\theta^2$ on $M$ satisfying
$$
\theta^1(T^1)\neq 0 \ \  
{\mbox{and}}\ 
\theta^2(T^2)\neq 0 .
$$
$$
\ker \theta^1_{\vert\ker\theta}=T^2 \  \  {\mbox{and}}\ 
 \  \ker \theta^2_{\vert\ker\theta}=T^1.
$$
Fixing one choice, all others are given by $\theta'^{i} = a^i \theta^i + v^i\theta$.

 On $E$ we define the tautological
form $\omega$.  That is $\omega_\theta=\pi^*(\theta)$ where
$\pi: {E}\rightarrow M$ is the natural projection.  
We also consider the tautological forms defined
by the forms $\theta^1$ and $\theta^2$ over the line bundle $E$.  That is, for each 
$\theta\in E$  we let $\omega^i_\theta= \pi^*(\theta^i)$.
At each point $\theta\in E$ we have the family of forms defined on $E$:
$$
        \omega' = \omega
$$
$$
        \omega'^{1} = a^1 \omega^1 + v^1\omega
$$
$$
        \omega'^{2} = a^2 \omega^2 + v^2\omega
$$

 We may, moreover, suppose that 
 $$
 d\theta= \theta^1\wedge \theta^2 \ \ \mbox{modulo} \ \theta
 $$ 
 and therefore
  $$
 d\omega= \omega^1\wedge \omega^2 \ \ \mbox{modulo} \ \omega .
 $$ 
  This imposes that $a^1a^2=1$.

  Those forms
vanish on vertical vectors, that is, vectors in the kernel
of the map $TE\rightarrow TM$. In order to define non-horizontal 1-forms 
we let $\theta$ be a section of $E$ over $M$ and introduce the coordinate
$\lambda\in \R^*$ in $E$.
By abuse of notation,
let $\theta$ denote the tautological form on the section $\theta$.
We write then the tautological form $\omega$ over $E$ is 
$$
\omega_{\lambda\theta}=\lambda \theta.
$$
Differentiating this formula we obtain
\begin{equation}  \label{domega}
                d\omega = 2\phi\wedge \omega + 
\omega^1\wedge\omega^{2}  
\end{equation}
where $\phi= \frac{d\lambda}{2\lambda}$ 
 modulo $\omega, \omega^1, 
\omega^{2}$.
Here
$
\frac{d\lambda}{2\lambda}
$ is a form intrinsically defined on $E$ up to horizontal
forms.
We obtain in this way a coframe bundle  satisfying equation \ref{domega} over $E$:
$$
        \omega' = \omega
$$
$$
        \omega'^{1} = a^1\omega^{1} + v^1 \omega
$$
$$
        \omega'^{2} = a^2\omega^{2} + v^2 \omega
$$
$$
        \phi'= \phi -\frac{1}{2}a^1v^{2}\omega^1 +\frac{1}{2}a^2v^1\omega^{2} +s\omega
$$
$v^1, v^2,s \in \R$  and $a^1,a^2\in \R^*$ such that $a^1a^2=1$.

\begin{dfn}
We denote by $Y$ the  coframe bundle
 $Y\rightarrow E$ given by the set of 1-forms  $\omega, \omega^{1},
  \omega^{2}, \phi$ as above.
  Two coframes
are related by
$$
(\omega', \omega'^{1},
  \omega'^{2}, \phi')=
(\omega, \omega^{1},
  \omega^{2}, \phi)
 \left ( \begin{array}{cccc}

                        1       &    v^1                &       v^2 & s \\

                        0 &     a^1 & 0            &       -\frac{1}{2}a^1 v^2  \\
                        0   &  0 & a^2             &       \frac{1}{2}a^2v^{ 1} \\
                        0      &  0 & 0     &       1

                \end{array} \right )
$$
where  and $s, v^1, v^2 \in \R$ and $a^1,a^2\in \R^*$ satisfy $a^1a^2=1$.
\end{dfn}

The bundle $Y$ can also be fibered over the manifold $M$.  In order to describe the bundle $Y$ as a principal fiber bundle over
 $M$ observe that choosing a local section $\theta$ of $E$ and forms $\theta^1$ and $\theta^2$ on $M$ such that
 $d\theta=\theta^1\wedge \theta^2$ one can write a trivialization of the fiber 

$$
      \omega=  \lambda \theta
$$
$$
      \omega^1=   a^1\theta^{1} + v^1\lambda\theta
$$
$$
       \omega^2= a^2\theta^{2} + v^2\lambda\theta
$$
$$
       \phi= \frac{d\lambda}{2\lambda} -\frac{1}{2}a^1v^{2}\theta^1 +\frac{1}{2}a^2v^1\theta^{2} +s\theta,
         $$
where $v^1, v^2,s \in \R$  and $a^1,a^2\in \R^*$ such that $a^1a^2=\lambda$.  Here the coframe $\omega, \omega^1,\omega^2,\phi$ is seen as 
composed of tautological forms.

The group $H$ acting on the right of this bundle 
 is 
$$
H=\left\{ \left ( \begin{array}{cccc}

                        \lambda   	&    	v^1 \lambda      	&       v^2\lambda 	& 		s \\

                       0				 &     a^1 			&		 0            &      - \frac{1}{2}a^1 v^2\\
                       0   				&  	0 				&		 a^2        &        \frac{1}{2}a^2v^{ 1}\\
                        0		       &  		0              &     0   			&       1

                \end{array} \right )
\mbox{
where   $s, v^1, v^2 \in \R$ and $a^1,a^2\in \R^*$ satisfy $a^1a^2=\lambda$
}
\right\}.
$$

Consider the homomorphism from the Borel group $B\subset \SL(3,\R)$ of upper triangular matrices  with determinant one into $H$
$$
j: B\rightarrow H
$$
given by

$$
\left( \begin{array}{ccc}

                        a  	&    	c      	&       e \\

                       0				 &     \frac{1}{ab}		&		f\\
                       0   				&  	0 				&		b
                      
       \end{array} \right ) 
       \longrightarrow  
       \left ( \begin{array}{cccc}

                        \frac{a}{b}  		&    	-a^2 f    	&       	\frac{c}{b}		& 		-eb+\frac 1 2acf \\

                       0				&     a^2b 		&		0            		&      		-\frac 1 2 abc\\
                       0   			&  	0 		&		\frac{1}{ab^2}    &        	-\frac{f}{2b}\\
                       0		         	&  	0              	&       	0   			&       	1

                \end{array} \right )
$$
One verifies that the homomorphism is surjective 
so  that $H$ is isomorphic to the Borel group of 
upper triangular matrices in $\SL(3,\R)$.

\begin{prop}  The bundle $Y\rightarrow M$ is a principal bundle with structure group $H$. 
\end{prop}\label{proposition:coframe}

\subsection{Construction of connection forms on the bundle $Y$}\label{section:connectionforms}

 The goal of this section is to  review the construction of canonical forms defined on the coframe bundle $Y\rightarrow E$ as in \cite{FV}.  They 
 give rise to a Cartan connection on $Y$  with values in $   {\mathfrak {sl}}(3,\R)$.

 A local section of the  coframe bundle over $M$ may be given by three forms 
 $$
 \theta, \ \ \theta^1, \ \ \theta^2
 $$
 satisfying 
 $d\theta=\theta^1\wedge \theta^2$, with 
$
\ker \theta^1_{\vert\ker\theta}=T^2 $ and  $\ker \theta^2_{\vert\ker\theta}=T^1.
$
They give coordinates on the  cotangent bundle over $E$.  
  Indeed, at  $\lambda\theta\in E$, the coframes of $Y$ are parametrized as follows:
 $$
 \omega= \lambda \theta
 $$
$$
 \omega^i=a^i\theta^i + v^i\lambda\theta
$$
 with $a^1a^2=\lambda$ and
  \begin{equation}\label{ome}
 d\omega =2 \phi\wedge \omega +\omega^1\wedge \omega^2,
  \end{equation}
 where $\phi = \frac{d\lambda}{2\lambda}  \  \mod{\omega^1,\omega^2,\omega}$.
 
 Differentiating the forms $\omega^1$ and $\omega^2$ we obtain new forms which correspond to the coordinates $a^1, v^1, a^2, v^2$ (recall that $a^1$ and $a^2$ are not independent):  

\begin{lem}
There exists linearly independent forms $w, \phi^1,\phi^2$ defined on $T^*Y$
such that
\begin{equation}\label{domegai}
d\omega^1=\phi \wedge  \omega^1+
3w \wedge \omega^1 +\omega\wedge \phi^1\  \mbox{and}\ \ 
d\omega^2=\phi \wedge  \omega^2 -
3 w\wedge \omega^2-\omega\wedge \phi^2
\end{equation}
with $w=\frac{1}{6}(\frac{da^1}{a^1}-\frac{da^2}{a^2})$ mod $(\omega,\omega^1,\omega^2)$ and $\phi^1=-dv^1$, $\phi^2=dv^2$ mod $(\omega,\omega^1,\omega^2)$.
\end{lem}

The coordinate $s$ in the bundle $Y$ is associated to a new form: 
\begin{lem}
There exists a 1-form $\psi$ such that
\begin{equation}\label{dphi}
d\phi = \omega \wedge \psi-\frac{1}{2}(
\phi^{2}\wedge \omega ^1+
\phi^1\wedge \omega^{2})
\end{equation}
\end{lem}

The forms $w, \phi^1,\phi^2$ and $\psi$ are not yet determined.  Define
$$
W=dw+\frac{1}{2} \omega^2\wedge \phi^1 -\frac{1}{2} \omega^1\wedge \phi^2
$$
$$
\Phi^1= d\phi^1 +3\phi^1\wedge w+\omega^1\wedge \psi+\phi\wedge \phi^1 
$$
$$
\Phi^2=d\phi^2 -3\phi^2\wedge w-\omega^2\wedge \psi+\phi\wedge \phi^2
$$ 
\begin{lem}
There exists unique forms $w, \phi^1, \phi^2$ and $\psi$ such that $W=0$, $\Phi^1=Q^1\omega\wedge \omega^2$ and  $\Phi^2=Q^2\omega\wedge \omega^1$ where $Q^1$ and $Q^2$ are functions on  $Y$.
\end{lem}
 
 We can represent the structure equations \ref{ome}, \ref{domegai}, \ref{dphi} 
  as a matrix equation whose entries are differential forms. The forms are disposed in the Lie algebra ${\mathfrak {sl}}(3,\R)$ and define a Cartan connection on $Y$.
  
  \begin{prop}\label{Cartan-path}
There exists a unique Cartan connection $\pi :TY\to {\mathfrak {sl}}(3,\R) $ defined on $Y$ of the form
$$\pi=\left ( \begin{array}{ccc}
                        \phi+w     &   \phi^{ 2}   &   \psi  \\

                     \omega^1          &   -2 w    &    \phi^1  \\

                      \omega         &  \omega^{2}   &      -\phi+w
                \end{array} \right ).
$$
such that its curvature satisfies
$$\Pi= d\pi+\pi\wedge\pi=\left ( \begin{array}{ccc}
                        0 &   \Phi^{2}      &  \Psi   \\

                         0  &  0    &  \Phi^1  \\

                           0  &   0  &    0
               \end{array} \right )
$$
with $
\Phi^1=Q^1\omega\wedge \omega^2,\ \ 
\Phi^2=Q^2\omega\wedge \omega^1 \ \mbox{and}\ \ 
\Psi =\left( U_1\omega^1+ U_2\omega^2\right) \wedge \omega.
$
\end{prop}

 \subsection{Curvature forms and Bianchi identities}

 Curvature forms appear as differentials of connection forms and are used implicitly in order to
 fix the connection forms.
 
   We recall:
  \begin{equation}\label{omega11}
W=dw-\frac{1}{2} \omega^2\wedge \phi^1 +\frac{1}{2} \omega^1\wedge \phi^2=0,
\end{equation}
 \begin{equation}\label{Phi1}
\Phi^1= d\phi^1 +3\phi^1\wedge w+\omega^1\wedge \psi+\phi\wedge \phi^1  =Q^1\omega\wedge \omega^2,
\end{equation}
 \begin{equation}\label{Phi2}
\Phi^2=d\phi^2 -3\phi^2\wedge w-\omega^2\wedge \psi+\phi\wedge \phi^2=Q^2\omega\wedge \omega^1,
\end{equation}
 \begin{equation}\label{dPsi'}
\Psi := d\psi-\phi^1\wedge \phi^2-2\phi\wedge \psi =(U_1\omega^1+ U_2\omega^2)\wedge \omega.
 \end{equation}
where  $Q^1, Q^2, U^1$ and $U^2$ are functions on $Y$.
  
  \subsubsection{}
  Equation $d(d\phi^1)=0$ obtained differentiating $\Phi^1$ above implies
   \begin{equation}\label{dq1}
   dQ^1 -6Q^1 w +4Q^1\phi=S^1\omega +U_2\omega^1 + T^1\omega^2,
   \end{equation}
`where we introduced functions $S^1$ and $T^1$.

   \subsubsection{}
   Analogously, equation $d(d\phi^2)=0$ obtained differentiating $\Phi^2$ above 
   implies
   \begin{equation}\label{dq2}
   dQ^2 +6Q^2 w +4Q^2\phi=S^2\omega -U_1\omega^2 + T^2\omega^1,
   \end{equation}
   where we introduced new functions $S^2$ and $T^2$.

\subsubsection{} Equation  $d(d\psi)=0$ obtained from 
\ref{dPsi'}
implies
 \begin{equation}\label{du1}
 dU_1+5U_1\phi +3U_1 w +Q^2\phi^1=A\omega+B\omega^1+C\omega^2
 \end{equation}
 and
 
  \begin{equation}\label{du2}
 dU_2+5U_2\phi -3U_2 w -Q^1\phi^2=D\omega+C\omega^1+E\omega^2.
 \end{equation}

 \subsection{Embedding $\iota_2 : Y_2\rightarrow Y$}\label{section:embedddingY2}

 
 The goal now is to obtain an immersion $\iota_2:Y_2\rightarrow Y$.  
 One can construct the bundle $Y_2$ using the bundle $E$ of contact forms as a first step.  Than $Y_2$ is a coframe bundle over $E$ obtained by the tautological forms $\omega, \omega^1, \omega^2$ corresponding to forms $\theta, \theta^1,\theta^2$ satisfying $d\theta=\theta^1\wedge \theta^2 + 2\phi\wedge \omega$ with an appropriate $\phi$.

 By abuse of language again as for $\iota_1: Y_1\rightarrow Y_2$, we may write the connection forms of $Y_1$ and $Y$ using the same letters and then obtain:
 
 \begin{prop}\label{propY2toY}
 There exists a unique embedding $\iota_2 : Y_2\rightarrow Y$ satisfying
 $\iota_2^*(\omega)=\omega$, $\iota_2^*(\omega^1)=\omega^1$, $\iota_2^*(\omega^2)=\omega^2$, $\iota_2^*(\phi)=\phi$.
 
 \end{prop}

 \Pf As $Y_2$ and $Y$ are both coframe bundles over the line bundle $E$ of all contact forms, we can assume that
 the embedding projects to the identity map on $E$.   Over $E$, $Y$ is a coframe bundle with structure group
 $$
\left\{ \left( \begin{array}{ccc}

                        a  	&    	c      	&       e \\

                       0				 &     \frac{1}{a^2}		&		f\\
                       0   				&  	0 				&		a
                      
       \end{array} \right ) \right\}.
 $$
 In order to determine the embedding we need to choose functions $c$, $e$ and $f$.  The diagonal matrix correspond to the fiber of $Y_2$ and does not need to be fixed.  Consider then a map from $M$ to the group above given by
 $$
h=  \left( \begin{array}{ccc}

                        1  	&    	c      	&       e \\

                       0				 &     {1}		&		f\\
                       0   				&  	0 				&		1
                      
       \end{array} \right ).
$$
Recall that
$$
{R_h}^*\pi= h^{-1}d\,h+ Ad_{h^{-1}}\pi.
$$
We obtain, neglecting the terms of the connection of $Y$ which are not relevant, the following transformation formulae.  Remark that the term $h^{-1}d\,h$ does not appear in the transformation of these components.

\begin{align}
\tilde\omega  &=  \omega\notag\\
\tilde{\omega}^1 &= \omega^1-f\,\omega\notag\\
\tilde\omega^2 &=\omega^2+{c}\,\omega\notag\\
\tilde\phi &=\phi-\frac 1 2 c\,\omega^1-{f}\,\omega^2+(\frac 1 2 cf-{e})\,\omega
\end{align}

The forms $\omega^1$ and $\omega^2$ defined at each point of $Y_2$ define corresponding forms 
$\omega^1$ and $\omega^2$ in $Y$.  We observe then that the functions  $f$ and $c$ must be zero in order that
$\iota_2^*(\omega^1)=\omega^1$, $\iota_2^*(\omega^2)=\omega^2$.  Finally the form $\phi$ on $Y_2$ defines a corresponding form on $Y$ and we conclude that $e=0$ if we impose that $\iota_2^*(\phi)=\phi$.

 \EPf
 
 \subsubsection{The curvature of $Y$ in terms of the curvature of $Y_2$}\label{section:embedddingY2-curvature}
 
 We obtain the following equations by pulling back to $Y_2$ the structure equations 
on $Y$ through the embedding $\iota_2$: 

$$
d\omega=2\phi \wedge \omega + \omega^1\wedge\omega^2
$$
\begin{equation}\label{dtheta1'}
d\omega^1=\phi \wedge  \omega^1+
3{\tilde{w}} \wedge \omega^1 +\omega\wedge \phi^1
\end{equation}
\begin{equation}\label{dtheta2'}
d\omega^2=\phi \wedge  \omega^2 -
3 {\tilde{w}}\wedge \omega^2-\omega\wedge \phi^2
\end{equation}
\begin{equation}
d\phi = \omega \wedge {\tilde{\psi}}-\frac{1}{2}(
\phi^{2}\wedge \omega ^1+
\phi^1\wedge \omega^{2})
\end{equation}
\begin{equation}\label{domega11'}
d{\tilde{w}}= -\frac{1}{2}\phi^2\wedge \omega^1+\frac{1}{2}\phi^1\wedge \omega^2
\end{equation}
\begin{equation}\label{phi1"}
d\phi^1 +3\phi^1\wedge {\tilde{w}}+\omega^1\wedge {\tilde{\psi}}+\phi\wedge \phi^1  =Q^1\omega\wedge \omega^2
\end{equation}
$$
d\phi^2 -3\phi^2\wedge {\tilde{w}}-\omega^2\wedge {\tilde{\psi}}+\phi\wedge \phi^2=Q^2\omega\wedge \omega^1
$$
$$
d{\tilde{\psi}}-\phi^1\wedge \phi^2-2\phi\wedge {\tilde{\psi}} =(U_1\omega^1+ U_2\omega^2)\wedge \omega.
$$

In the formulae above we write the pull back of any form $\alpha$ defined on $Y$ using the same notation $\alpha$ except for the pull backs $\tilde{w}=\iota_2^*w$ and ${\tilde{\psi}}=\iota_2^*\psi$.  
We should compare with the structure equations of $Y_2$ and obtain an expression for $Q^1$ and $Q^2$:
$$
d\omega=2\phi \wedge \omega + \omega^1\wedge\omega^2
$$
$$
d\omega^1=\phi \wedge  \omega^1+3w\wedge \omega^1 +\omega \wedge \tau^1
$$
$$
d\omega^2=\phi\wedge \omega^2-3w\wedge \omega^2 -\omega\wedge \tau^2
$$
with  $\tau^1\wedge \omega^2=\tau^2\wedge \omega^1=0$.

Recall also that $d\phi = 
\omega \wedge \psi$ with $\psi=A\omega^1+B\omega^2$ and $
dw=C\omega\wedge\omega^1+ D\omega\wedge\omega^2+S\omega^1\wedge\omega^2$, where $A,B,C,D$ and $S$ are functions on $Y_2$.  

\begin{itemize}

\item The differences between the pull back equations and the structure equations for $d\omega^1$ and $d\omega^2$ give, respectively,
$$
3({\tilde{w}-w} )\wedge \omega^1 +\omega\wedge( \phi^1-\tau^1)=0
$$
and 
$$
3({\tilde{w}-w} )\wedge \omega^2 +\omega\wedge( \phi^2-\tau^2)=0.
$$
Therefore, by Cartan's lemma 
$${\tilde{w}-w}=m\omega$$ for a function $m$ to be determined and
$$
\phi^1-\tau^1=-3m\omega^1+n\omega \ \mbox{and} \ \phi^2-\tau^2=-3m\omega^2+P\omega,
$$
where $n$ and $P$ are functions to be determined.
\item
The difference between the equation  $d\phi = 
\omega \wedge \psi$ and the pull back equation for $d\phi$ above is
$$
 \omega \wedge ({\tilde{\psi}}-\psi)-\frac{1}{2}(
\phi^{2}\wedge \omega ^1+
\phi^1\wedge \omega^{2})=0.
$$
Substituting the expressions for $\phi^1$ and $\phi^2$ obtained in the item above we obtain
$$
 \omega \wedge ({\tilde{\psi}}-A\omega^1-B\omega^2)-\frac{1}{2}(
(\tau^2-3m\omega^2+P\omega)\wedge \omega ^1+
(\tau^1-3m\omega^1+n\omega)\wedge \omega^{2})=0,
$$
which simplifies to
$$
\omega\wedge({\tilde{\psi}}-A\omega^1-B\omega^2-\frac{P}{2}\omega^1-\frac{n}{2}\omega^2)=0.
$$
This implies that
$$
{\tilde{\psi}}=(A+\frac{P}{2})\omega^1+(B+\frac{n}{2})\omega^2+q\omega,
$$
where $q$ is a function to be determined.

\item  The difference between the equations for $d\tilde{w}$ and $dw$ is then
$$
d(m\omega)=d\tilde{w}-dw=-\frac{1}{2}\phi^2\wedge \omega^1+\frac{1}{2}\phi^1\wedge \omega^2-S\omega^1\wedge \omega^2-C\omega\wedge \omega^1-D\omega\wedge \omega^2.
$$
Substituting in this formula the expressions for $\phi^1$ and $\phi^2$ in terms of the enriched structure we obtain
$$
dm\wedge \omega+m(2\phi\wedge\omega +\omega^1\wedge \omega^2)=(-S-3m)\omega^1\wedge \omega^2+(-\frac{P}{2}-C)\omega\wedge \omega^1+(\frac{n}{2}-D)\omega\wedge \omega^2.
$$
That is,
$$
(S+4m)\omega^1\wedge \omega^2+\omega\wedge \left(-dm-2m\phi+(\frac{P}{2}+C)\omega^1-(\frac{n}{2}-D)\omega^2\right)=0.
$$
Therefore $$m=-S/4$$ and so
$$
\frac{1}{4}dS+\frac{1}{2}S\phi+(\frac{P}{2}+C)\omega^1-(\frac{n}{2}-D)\omega^2+E\omega=0
$$
where $E$ is a function determined by the derivative of $S$.
Writing $dS+2S\phi=S_0\omega+S_1\omega^1+S_2\omega^2$ and comparing with the above expression, we obtain 
\begin{equation}\label{S0S1S2}
S_0=-4E, \ \ S_1=-2(P+2C), \ \ S_2=-2(-n+2D).
\end{equation}
 Therefore the functions  $P$ and $n$ are determined.
It remains to determine the function $q$.

\item Computing $d\phi^1= d(\frac{3S}{4}\omega^1+n\omega+\tau^1)$ and equating to the structure equation
$d\phi^1= -3\phi^1\wedge {\tilde{w}}-\omega^1\wedge {\tilde{\psi}}-\phi\wedge \phi^1  +Q^1\omega\wedge \omega^2$ we obtain, after a computation
writing 
\begin{equation}\label{dn}
dn+3n(\phi-w)=n_0\omega+n_1\omega^1+n_2\omega^2,
\end{equation}
\begin{equation}\label{n2}
n_1=-\tau^1_2\tau^2_1-q+\frac{9}{16}S^2-3E \ \ ,
n_2=-Q^1+\frac{3}{2}S\tau^1_2+\tau^1_{20}.
\end{equation}
Recalling that $n= S_2/2+4D$, $S$ and $D$ are determined by $Y_2$, we obtained an expression for 
$Q^1$ in terms of the enriched structure.  Note also that $q$  is determined by the first equation.

\item Analogously, computing $d\phi^2=-3\frac{3S}{4}\omega^2+p\omega +\tau^2$ and equating to the structure equation
$d\phi^2 -3\phi^2\wedge {\tilde{w}}-\omega^2\wedge {\tilde{\psi}}+\phi\wedge \phi^2=Q^2\omega\wedge \omega^1$ we obtain, after a computation, writing $dP+3P(\phi+w)=P_0\omega+P_1\omega^1+P_2\omega^2$),
\begin{equation}\label{p1}
P_1=-Q^2 -\frac{3}{2}S\tau^2_1+\tau^2_{10}\ \ ,
P_2=\tau^1_2\tau^2_1+q-\frac{9}{16}S^2-3E.
\end{equation}
Recalling that $P= -S_1/2-2C$, $S$ and $C$ are determined by $Y_2$, we obtained an expression for 
$Q^2$ in terms of the enriched structure.

\end{itemize}

The following proposition follows directly from the computations above.

\begin{prop} Suppose $Y_2$ with its adapted  Cartan connection has null torsion, that is,  satisfies $\tau^1=\tau^2=0$ and 
$dw=S\omega^1\wedge\omega^2$.
Then 
$$
Q^2=\frac{1}{2}S_{11}
$$
and
$$
Q^1=-\frac{1}{2}S_{22},
$$
where $S_{11}$ is the $\omega^1$ component of the form $dS_1$ and $S_{22}$ is the $\omega^2$ component of the form $dS_2$. 

\end{prop}

\Pf
From Proposition \ref{prop:nultorsion}, null torsion and the condition  that $dw=S\omega^1\wedge \omega^2$ (that is, $C=D=0$) implies that
$P=-S_1/2$ and $n= S_2/2$.  The result is therefore implied from the previous formulas.
\EPf

\subsection{The embedding $Y_1\to Y$}\label{section:embedding1}
 
 Recall that $Y_1$ is a coframe bundle of forms $(\theta,\theta^1,\theta^2)$  over $M$. 
 Choosing a local section, the pullback forms over $M$ are also denoted by $(\theta,\theta^1,\theta^2)$. We recall here an embedding $Y_1\rightarrow Y$ obtained in \cite{FV}.  
 
 A section $(\theta,\theta^1,\theta^2)$ of the coframe bundle $Y_1$ clearly defines a path geometry on $M$. We obtain then a  line bundle  $E$ and a principal bundle $Y$ with its associated Cartan connection.  Also, $(\theta,\theta^1,\theta^2)$ defines, up to the action by the group of matrices
 $$
 \left( \begin{array}{ccc}

                        1  	&    	c      	&       e \\

                       0				 &     1		&	f	\\
                       0   				&  	0 				&		1
                      
       \end{array} \right )
 $$
 sections of the tautological forms $(\omega,\omega^1,\omega^2)$ on $Y$.  In order to define a canonical section we use the following
 
 \begin{prop}
Let $\theta,\theta^1,\theta^2$ be a coframe section of $Y_1$ and consider the principal bundle $Y$ defined by this coframe.  Then
there exists a  unique section $s:M\to Y$ such that $s^* \omega=\theta$, $s^* \omega^1=\theta^1$, $s^* \omega^2=\theta^2$ and $s^*\phi=0$. 
 \end{prop}
 
It is easy to verify that this definition is equivariant with respect to  the action $G_1$, the one parameter group of the strict contact structure.
 This  defines then the embedding $Y_1\rightarrow Y$.

 \subsection{The equivalence problem for a second order differential equation}

In this section we recall the treatment by Cartan of the point equivalence between second order differential equations.  It is included in order to fix conventions and to compare the 
invariants defined in the next section.

Recall that for a second order differential equation we define
$$
\theta= d y-pdx,
$$
and
$$
L_1= \ker\{d p-Fdx\}\cap \ker\{dy-pdx\},\ \ \ L_2= \ker dx\cap \ker dy.
$$

For $Z^1=dx$ and $Z^2=d p-Fdx$, one has then $d\theta= Z^1\wedge Z^2$.
The general forms defining the lines at each tangent space may be described  by 
$$
\omega^1=a_1Z^1, \ \ \omega^2=a_2Z^2, \\ \omega=a_1a_2\theta
$$
where $a_1,a_2$ are non-vanishing positive functions on the manifold,
so that we have always 
$$2\phi\wedge\omega+\omega^1\wedge\omega^2=d\omega=(\frac{da_1}{a_1}+\frac{da_2}{a_2})\wedge\omega+a_1Z^1\wedge a_2Z^2,$$
and we obtain comparing with \ref{do} 
$$
\phi=\frac 1 2(\frac{da_1}{a_1}+\frac{da_2}{a_2})+r\omega.
$$

One computes 
\begin{equation}
(\phi+3w)\wedge\omega^1+\omega\wedge\tau^1=d\omega^1=da_1\wedge Z^1+a_1.0=\left(\frac 1 2(\frac{da_1}{a_1}+\frac{da_2}{a_2})+\frac1 2(\frac{da_1}{a_1}-\frac{da_2}{a_2})\right)\wedge \omega^1 
\end{equation}
and obtain
$$
3w=\frac1 2(\frac{da_1}{a_1}-\frac{da_2}{a_2})-r\omega+s\omega^1, \\ \tau^1_2=0.
$$

Observe that, if 
$f(x,y,p)$ then 
$$df=\frac{1}{a_1}\frac{df}{dx}\omega^1+\frac{1}{a_2}f_p\omega^2+  \frac{1}{a_1a_2}f_y\omega,$$ 
where $\frac{df}{dx} =f_x+f_yp+f_pF $.                                   
Also
\begin{equation}
\begin{array}{rcl}
d\omega^2&=&\frac{da_2}{a_2}\wedge \omega^2+a_2(-\frac{1}{a_2}F_p\omega^2-  \frac{1}{a_1a_2}F _y\omega)\wedge\frac{1}{a_1}\omega^1\\
&=&(\frac 1 2(\frac{da_1}{a_1}+\frac{da_2}{a_2})-\frac1 2(\frac{da_1}{a_1}-\frac{da_2}{a_2}))\wedge\omega^2+\frac{1}{a_1}F_p\omega^1\wedge\omega^2-  \frac{1}{(a_1)^2}F _y\omega\wedge\omega^1\\
&=&(\phi-3w-2r\omega+(s+ \frac{1}{a_1}F _p)\omega^1)\wedge\omega^2-\omega\wedge \frac{F_y}{(a_1)^2}\omega^1
\end{array}
\end{equation}
Comparing with $d\omega^2=(\phi-3w)\wedge\omega^2-\omega^\wedge\tau^2$ we obtain $r=0,s=-\frac{F_p}{a_1},\tau^2_1=\frac{F_y}{a_1^2}$,
$$
\phi=\frac 1 2(\frac{da_1}{a_1}+\frac{da_2}{a_2}), \\  3w=\frac1 2(\frac{da_1}{a_1}-\frac{da_2}{a_2})-\frac{F_p}{a_1}\omega^1.
$$
From $d\phi=0$, we obtain $\psi=0$, or $A=B=0$, and from $\tau^1_2=0$ it follows $D=\tau^1_{20}=\tau^1_{22}=0$.

From above we get
$$3dw=-dF_p\wedge \frac{1}{a_1}\omega^1=-(\frac{F_{pp}}{a_2}\omega^2+\frac{F_{py}}{a_1a_2}\omega)\wedge\frac{1}{a_1}\omega^1$$
and comparing with $dw=C\omega\wedge\omega^1+S\omega^1 \wedge\omega^2$ we obtain
$$3C=-\frac{F_{py}}{a_1^2a_2}, \ \ \  3S=\frac{F_{pp}}{a_1a_2}.$$
Also
$$
d\tau^2_1=d(\frac{F_y}{a_1^2})=-\frac{2}{a_1^2}(\phi+3w+\frac{F_p}{a_1}\omega^1)F_y+\frac{1}{a_1^2}(\frac{1}{a_1}\frac{dF_y}{dx}\omega^1+\frac{1}{a_1a_2}F_{yy}\omega+\frac{1}{a_2}F_{yp}\omega^2).
$$
Comparing with $d\tau^2_1=-2\tau^2_1(\phi+3w)-3C\omega^2+\tau^2_{10}\omega+\tau^2_{11}\omega^1$ we obtain
$$
\tau^2_{11}=\frac{1}{a_1^3}(-2F_pF_y+\frac{dF_y}{dx}), \ \ \  \tau^2_{10}=\frac{1}{a_1^3a_2}F_{yy}.
$$
Now
$$
3dS=-2(3S\phi)+\frac{1}{a_1a_2}(\frac{1}{a_1}\frac{dF_{pp}}{dx}\omega^1+\frac{1}{a_2}F_{ppp}\omega^2+\frac{1}{a_1a_2}F_{ppy}\omega)
$$
and comparing with $dS=-2S\phi+S_0\omega+S_1\omega^1+S_2\omega^2$ we obtain
$$S_1=\frac{1}{3a_1^2a_2}\frac{dF_{pp}}{dx}, \ \ \ S_2=\frac{1}{3a_1a_2^2}F_{ppp}, \\ S_0=\frac{1}{3a_1^2a_2^2}F_{ppy}.$$
 It follows from $S_1=-2P-4C$  that $6P=\frac{1}{a_1^2a_2}(4F_{yp}-\frac{dF_{pp}}{dx}).$ Then
$$6dP=-\frac{1}{a_1^2a_2} (3\phi+3w+\frac{F_p}{a_1}\omega^1)(4F_{yp}-\frac{dF_{pp}}{dx})+
\frac{4}{a_1^2a_2}(\frac{1}{a_1}\frac{dF_{yp}}{dx}\omega^1+\frac{1}{a_2}F_{ypp}\omega^2+\frac{1}{a_1a_2}F_{yyp}\omega)$$
$$-\frac{1}{a_1^2a_2}(\frac{1}{a_1}\frac{d^2F_{pp}}{dx^2}\omega^1+\frac{1}{a_2}(\frac{dF_{ppp}}{dx}+F_{ypp}+F_{ppp}F_p)\omega^2+\frac{1}{a_1a_2}(\frac{dF_{ypp}}{dx}+F_{ppp}F_y)\omega)
$$
Comparing with $dP=-(3\phi+3w)P+P_0\omega+P_1\omega^1+P_2\omega^2$ we obtain
$$
P_0=\frac{1}{6a_1^3a_2^2}(4F_{yyp}-\frac{dF_{ypp}}{dx}-F_{ppp}F_y),  \ \ \ \ P_1= \frac{1}{6a_1^3a_2}(-4F_pF_{yp}+F_p\frac{dF_{pp}}{dx}+4\frac{dF_{yp}}{dx}-\frac{d^2F_{pp}}{dx^2}),
$$			
$$
P_2=\frac{1}{6a_1^2a_2^2}(4F_{ypp}-\frac{dF_{ppp}}{dx}-F_{ypp}-F_{ppp}F_p).
$$
From $Q^2=\tau^2_{10}-\frac 3 2 S\tau^2_1-P_1$ it follows
\begin{equation}\label{q2}
Q^2=\frac{1}{6a_1^3a_2}(6F_{yy}-3F_yF_{pp}+4F_pF_{yp}-F_p\frac{dF_{pp}}{dx}-4\frac{dF_{yp}}{dx}+\frac{d^2F_{pp}}{dx^2}).
\end{equation}
It follows from $S_2=2n-4D$ that $6n=\frac{1}{a_1a_2^2}F_{ppp}$. Then
$$6dn=-\frac{1}{a_1a_2^2}(3\phi-3w-\frac{F_p}{a_1}\omega^1)F_{ppp}+\frac{1}{a_1a_2^2}(\frac{1}{a_1}\frac{dF_{ppp}}{dx}\omega^1+\frac{1}{a_2}F_{pppp}\omega^2+\frac{1}{a_1a_2}F_{pppy}\omega).
$$
Comparing with $dn=-n(3\phi-3w)+n_0\omega+n_1\omega^1+n_2\omega^2$ we obtain
$$
n_0=\frac{1}{6a_1^2a_2^3}F_{pppy}, \ \ \ n_1=\frac{1}{6a_1^2a_2^2}(\frac{dF_{ppp}}{dx}+F_pF_{ppp})\ \ \ n_2=\frac{1}{6a_1a_2^3}F_{pppp}.
$$
From $Q^1=\tau^1_{20}+\frac{3}{2}S\tau^1_2-n_2$ it follows
\begin{equation}\label{q1}
Q^1=-\frac{1}{6a_1a_2^3}F_{pppp}.
\end{equation}
Formulas \ref{q1} and \ref{q2} are in \cite{Car}.

\section{A global invariant}\label{section:global}

In this section we define the global invariant for path structures.  It has a very similar definition with  the global invariant obtained in \cite{FV} in the context of a structure defined on the complexified tangent space of a 3-manifold.  But we  make the definition explicit in the case of path structures for the sake of clarity and to adapt differences of conventions with our previous paper.

Define the second Chern class of the bundle $Y$ with connection form $\pi$ as
$$
c_2(Y,\pi)=\frac{1}{8\pi^2}\tr(\Pi\wedge \Pi).
$$

$$
\left ( \begin{array}{ccc}
                        0 &   \Phi^{2}      &  \Psi   \\

                         0  &  0    &  \Phi^1  \\

                           0  &   0  &    0
               \end{array} \right )
               \wedge 
               \left ( \begin{array}{ccc}
                        0 &   \Phi^{2}      &  \Psi   \\

                         0  &  0    &  \Phi^1  \\

                           0  &   0  &    0
               \end{array} \right )
               =\left ( \begin{array}{ccc}
                        0 &   0     &  \Phi^1\wedge \Phi^2  \\

                         0  &  0    &  0\\

                           0  &   0  &    0
               \end{array} \right ).
$$
As $\Phi^1=Q^1\omega\wedge \omega^2$ and $\Phi^2=Q^2\omega\wedge \omega^1$ we have $\Pi\wedge \Pi=0$ and therefore
$$
c_2(Y,\pi)=0.
$$

\begin{dfn}
The transgression form is defined as
$$
TC_2(\pi)=\frac{1}{8\pi^2}\left ( \tr(\pi\wedge \Pi)+\frac{1}{3}\tr(\pi\wedge \pi\wedge\pi)\right)=
\frac{1}{24\pi^2}\tr(\pi\wedge \pi\wedge\pi).
$$
\end{dfn}

\begin{lem} The transgression form is closed, that is,  $d\, TC_2(\pi)= c_2(Y,\pi)=0$.

\end{lem}
\Pf 
We compute first, using the expressions of $\Phi^1$, $\Phi^2$ and $\Psi$, that
 $$
\tr (\Pi\wedge \pi) =\Phi^2\wedge \omega^1+\Phi^1\wedge \omega^2+\Psi\wedge\omega=0.
$$
Differentiating the curvature form we obtain $d\,\Pi=\Pi\wedge \pi-\pi\wedge \Pi$ and therefore

$$
0=d\,\tr (\Pi\wedge \pi)=\tr(d\,\Pi\wedge \pi+\Pi\wedge d\,\pi)
=\tr((\Pi\wedge \pi-\pi\wedge\Pi)\wedge \pi+\Pi\wedge (\Pi-\pi\wedge\pi))
$$
$$
=-\tr(\pi\wedge\Pi\wedge \pi).
$$

Note that $\tr(\alpha\wedge \beta)=(-1)^{kl}\tr(\beta\wedge \alpha)$ if $\alpha$ and $\beta$ are two matrices of forms 
of degree $k$ and $l$ respectively.  Therefore, computing
$$
\frac{1}{3}d\,\tr(\pi\wedge \pi\wedge\pi)= \tr(d\,\pi\wedge \pi\wedge\pi)= \tr((\Pi-\pi\wedge\pi)\wedge \pi\wedge\pi)
$$
$$
=- \tr(\pi\wedge\pi\wedge \pi\wedge\pi)=0.
$$
\EPf

\begin{dfn}\label{definition:global}
Suppose  that the fiber bundle $Y\rightarrow M$ is trivial and let $s: M\rightarrow Y$ be a section, we define then
$$
\mu=\int_M s^* TC_2(\pi)= \frac{1}{24\pi^2}\int_M s^*\tr(\pi\wedge \pi\wedge\pi).
$$
\end{dfn}

In principle that integral depends on the section but the following proposition shows that the integrand 
$$
s^* TC_2(\pi)
$$
defines
an element in $H^3(M,\R)$ which does not depend on the section.

  \begin{prop} 
 Suppose $s$ and $\tilde s$ are two sections.  Then
$$
 \tilde{ s}^*TC_2(\pi)-s^*TC_2(\pi)=-\frac{1}{8\pi^2}d\, s^*\tr (h^{-1}\pi\wedge d\,h ).
$$
where $h: M\rightarrow H$ is a map such that
 $\tilde s = R_h\circ s$.  
 \end{prop}
 \Pf
 Fix the section $s$. Than there exists a map $h: M\rightarrow H$ such that
 $\tilde s = R_h\circ s$.  We have then
 $$
 \tilde s^*TC_2(\pi)=\frac{1}{24\pi^2}s^*\tr(R_h^*\pi\wedge R_h^*\pi\wedge R_h^*\pi).
 $$
From the formula
$$
{R_h}^*\pi= h^{-1}d\,h+ Ad_{h^{-1}}\pi,
$$
we obtain
$$
\tr(R_h^*\pi\wedge R_h^*\pi\wedge R_h^*\pi)=
$$
$$
\tr \left( h^{-1}d\,h\wedge  h^{-1}d\,h\wedge h^{-1}d\,h 
+ 3h^{-1}d\,h\wedge  h^{-1}\pi\wedge d\,h+3h^{-1}\pi\wedge \pi\wedge d\,h + \pi\wedge \pi\wedge\pi \right)
$$
$$
=\tr \left( - h^{-1}d\,h\wedge  d\,h^{-1}\wedge d\,h 
-3d\,h^{-1}\wedge \pi\wedge d\,h+3h^{-1}\pi\wedge \pi\wedge d\,h + \pi\wedge \pi\wedge\pi \right).
$$

 Observe that the first term in the right hand side vanishes.  Indeed, $ d\,h^{-1}\wedge d\,h$ is upper triangular with null diagonal.  Moreover $h^{-1}d\,h$ is upper triangular and therefore the Lie algebra valued form also has zero diagonal.
Therefore 
$$\tr (h^{-1}d\,h\wedge  d\,h^{-1}\wedge d\,h )=0.
$$
  By the same argument $ \tr \left( h^{-1}\Pi \wedge d\,h\right)=0$.

Now we show that  $$d\,\tr (h^{-1}\pi\wedge d\,h )=
\tr \left( d\,h^{-1}\wedge \pi\wedge d\,h-h^{-1}\pi\wedge \pi\wedge d\,h\right).$$

Compute
$d\tr (h^{-1}\pi\wedge d\,h )=\tr \left( d\,h^{-1}\wedge \pi\wedge d\,h+h^{-1}d\pi\wedge d\,h\right)$
$$
=\tr \left( d\,h^{-1}\wedge \pi\wedge d\,h+h^{-1}(\Pi -\pi\wedge \pi)\wedge d\,h\right),
$$
which gives, using that $ \tr \left( h^{-1}\Pi \wedge d\,h\right)=0$,
$$
d\tr (h^{-1}\pi\wedge d\,h )=\tr \left( d\,h^{-1}\wedge \pi\wedge d\,h-h^{-1}\pi\wedge \pi\wedge d\,h\right).
$$
We obtained therefore that
$$
 \tilde{e} s^*TC_2(\pi)=s^*TC_2(\pi)-\frac{1}{8\pi^2}d\, s^*\tr (h^{-1}\pi\wedge d\,h )
$$
and this completes the proof of the proposition.

\EPf

Let $\mu(t)$ be the invariant defined as a function of the a parameter describing the deformation of the structure 
on a closed manifold $M$
and define $\delta\mu=\frac{d}{dt}\mu (0)$.  One can interpret the flat structures as giving critical points of 
the global invariant $\mu$ through the first variation formula which we refer to \cite{FV} for a proof.

\begin{prop}\label{proposition:variation}
$$
\delta\mu= -\frac{1}{4\pi^2}\int_M s^*\tr(\dot \pi\wedge \Pi).
 $$
 \end{prop}


 {The global invariant can be computed most easily for a path structure induced by an enriched or strict path structure.}
 
 \begin{prop}\label{proposition-invariantY_2}
 Let $M$ be an enriched path structure and $Y_2\rightarrow Y$ be the canonical embedding of the enriched geometry into the induced path geometry.  Then 
 $$
8\pi^2\mu(Y)=
(n_1-\frac 3 4S_0+2\tau^1_2\tau^2_1)\omega\wedge\omega^1\wedge\omega^2+\omega\wedge\omega^1\wedge(-2A\phi-(\frac 3 2S_1+6C)w)
$$
$$
+\omega\wedge\omega^2(-2B\phi-(\frac 3 2 S_2+6D)w)-\frac 9 2 Sw\wedge\omega^1\wedge\omega^2.
$$
 \end{prop}
 \Pf
 One compute first the following formula.
 $$
\frac 1 3 \tr (\pi\wedge\pi\wedge \pi)=(2\omega\wedge \phi -\omega^1\wedge \omega^2)\wedge \psi
-\omega\wedge \phi^1\wedge \phi^2+\omega^1\wedge (\phi+3w)\wedge\phi^2
+\omega^2\wedge (\phi-3w)\wedge\phi^1.
$$
Therefore using the embedding of $Y_2\rightarrow Y$ in the previous section we obtain by a computation:
$$
\frac 1 3 \tr (\pi\wedge\pi\wedge \pi)=(n_1-\frac 3 4S_0+2\tau^1_2\tau^2_1)\omega\wedge\omega^1\wedge\omega^2+\omega\wedge\omega^1\wedge(-2A\phi-(\frac 3 2S_1+6C)w)
$$
$$
+\omega\wedge\omega^2(-2B\phi-(\frac 3 2 S_2+6D)w)-\frac 9 2 Sw\wedge\omega^1\wedge\omega^2.
$$
\EPf

Using the embedding of $Y_1\rightarrow Y$ in the previous section we obtain by a similar computation:
\begin{prop}
 Let $M$ be a strict  path structure and $Y_1\rightarrow Y$ be the canonical embedding of the strict geometry into the induced path geometry.  Then 
\begin{equation}\label{pi3}
8\pi^2\mu(Y)=
(n_1-\frac 3 4S_0+2\tau^1_2\tau^2_1)\omega\wedge\omega^1\wedge\omega^2-w\wedge((\frac 3 2S_1+6C)\omega\wedge\omega^1+(\frac 3 2 S_2+6D)\omega\wedge\omega^2+\frac 9 2 S\omega^1\wedge\omega^2)
\end{equation}
\end{prop}

\subsection{The global invariant for second order differential  equations on the torus}\label{section:globalDE}

In this section we obtain formulas for the global invariant in the case of an ordinary differential equation defined on the torus.
Recall that the projectivized cotangent bundle $\pi : PT^*S\to S$ of a surface $S$ is described locally by $(x,y,[p,q])$ where $(x,y)$ are local coordinates on the surface 
and $pdx+qdy$ is a form at $(x,y)$.  The Liouville form $\theta$ on $T^*S$ is defined to be the tautological form $\theta (x,y,-pdx+qdy)=\pi^*(-pdx+qdy)$.  It induces
a contact distribution on $PT^*S$, which in the chart $(x,y,p)\to (x,y,[p,1])$ is given by the kernel of the form $dy-pdx$.  On the chart $(x,y,q)\to (x,y,[1,q])$ the contact distribution is the kernel of $dx-qdx$.
One can also consider, fixing a metric on the surface, the unit cotangent bundle $(T^*)^1S$ which is a double cover of $PT^*S$.

The fibers of the bundle $PT^*S$ give a canonical field of directions on the Liouville distribution.  Observe that,  in local coordinates $(x,y,p)$, it is described by
$\ker dx\cap \ker dy= \ker dx\cap \ker(dy-pdx)$.  Choosing another direction on the contact distribution amounts to define
a form, in local coordiantes $(x,y,p)$,  $dp-G(x,y,p)dx$, where $G(x,y,p)$ is a function.  On the chart $(x,y,q)$ one writes then
$$
d(\frac{1}{q})-G(x,y,\frac{1}{q})dx=-\frac{1}{q^2}dq-G(x,y,\frac{1}{q})dx.
$$
Therefore, the direction is determined by $dq+G(x,y,\frac{1}{q})q^3dy$ (the contact distribution  is $\ker(dx-qdy)$).  In order to have a well defined 
direction we need that the function $G(x,y,\frac{1}{q})q^3$ has a differentiable extension for $q=0$.

\begin{dfn}
A second order differential equation on a surface $S$ is a path structure on the projective cotangent bundle with contact structure induced by the Liouville form
and such that one of the directions is given by the fibers. 
\end{dfn}

It is convenient to introduce a new coordinate in the fiber $\alpha \in ]-\pi,\pi]$ through the  formula $p=\tan{\alpha/2}$.  The contact distribution is defined by a globally defined form on the coordinates $(x,y,\alpha)$:
$$
\theta=\cos{\alpha/2}dy-\sin{\alpha/2}dx.
$$
The fiber direction is defined by the equations $dx=dy=0$ which can also be described by, defining 
$\theta^1= \sin{\alpha/2}dy+\cos{\alpha/2}dx$, as $\ker \theta^1\cap \ker \theta$.  
The last form, which depends on a choice of  a function,  is
$$
\theta^2=d\alpha -F(x,y,\alpha)\theta^1.
$$
Observe that 
$$
d\theta= \frac{1}{2} \theta^1\wedge \theta^2.
$$
The relation with the differential equation given on the chart $(x,y,p)$ is given writing
$$
dp=\frac{1}{2}(1+p^2)d\alpha
$$
and therefore, as $dy=pdx$ in that chart,
$$
d\alpha-F(x,y,\alpha)\theta^1= \frac{2}{1+p^2}dp-F(x,y,2\arctan{p})(\sin{\alpha/2}dy+\cos{\alpha/2}dx)
$$
$$
= \frac{2}{1+p^2}dp-F(x,y,2\arctan{p})(\sin{\alpha/2}. pdx+\cos{\alpha/2}dx)
$$
$$
=\frac{2}{1+p^2}dp-F(x,y,2\arctan{p})(\sin{\alpha/2}.\tan{\alpha/2}+\cos{\alpha/2})dx.
$$
and recalling that $\cos{\alpha/2}=\frac{1}{\sqrt{1+p^2}}$,

$$
=\frac{2}{1+p^2}dp-F(x,y,2\arctan{p}){(1+p^2)^{1/2}}dx.
$$
Therefore $2G(x,y,p)= F(x,y,2\arctan{p}){(1+p^2)^{3/2}}$.

\vspace{1cm}

\subsubsection{The strict and enriched structure of a differential equation on the torus}

Here we will work with a double cover of the projective cotangent bundle of the torus.
We define the path structure associated to a differential equation on the torus 
through a strict path structure defined by
$$
\theta=\cos{\alpha}dy-\sin{\alpha}dx.
$$
$$
\theta^1= \sin{\alpha}dy+\cos{\alpha}dx
$$ 
and
$$
\theta^2=d\alpha -F(x,y,\alpha)\theta^1.
$$
Here $F(x,y,\alpha)$ is a function defined on the torus.  Observe that 
$$
d\theta=  \theta^1\wedge \theta^2.
$$

In the following we will write, for a function $f: T^3\to \R$, defined on the torus,
$$
df=f_0\theta+ f_1\theta^1+f_2\theta^2,
$$
so that $f_x=-f_0\sin{\alpha}+(f_1-f_2 F)\cos{\alpha}$, $f_y=f_0\cos{\alpha}+(f_1-f_2 F)\sin{\alpha}$ and
$f_\alpha=f_2$.
Compute 
$$
d\theta^1
=(\theta^2+F\theta^1)\wedge \theta
$$
and
$$
d\theta^2
=((F_0-F^2)\theta^1-F\theta^2)\wedge \theta +F_2\theta^1\wedge \theta^2.
$$

Consider now the enriched structure defined by $\theta, \theta^1$ and $\theta^2$ and the tautological forms
 $\omega=a_1a_2\theta, \omega^1=a_1\theta^1$ and $\omega^2=a_2\theta^2$.
 
 We first compute
 $$
 d\omega=2\phi\wedge \omega+\omega^1\wedge \omega^2,
 $$
where $\phi=\frac 1 2(\frac{da_1}{a_1}+\frac{da_2}{a_2})$.
Next we compute $d\omega^1$ and $d\omega^2$:
$$
d\omega^1=\frac{da_1}{a_1}\wedge a_1\theta^1+ a_1(\theta^2+F\theta^1)\wedge \theta
$$
$$
d\omega^2=\frac{da_2}{a_2}\wedge a_2\theta^2+ a_2(((F_0-F^2)\theta^1-F\theta^2)\wedge \theta +F_2\theta^1\wedge \theta^2)
$$
Comparing with the structure equations of the enriched structure in \ref{domegaitau} we may write
\begin{equation}
d\omega^1=\phi \wedge \omega^1+
3 w\wedge \omega^1 +\omega \wedge \tau^1\ \ \ 
d\omega^2=\phi \wedge \omega^2-
3w\wedge \omega^2-\omega\wedge \tau^2.
\end{equation}
with
$$
3w= \frac 1 2(\frac{da_1}{a_1}-\frac{da_2}{a_2})- \frac{1}{a_1a_2}  F\omega- \frac{1}{a_1}  F_2\omega^1,
$$
$$
\tau^1=-\frac{1}{a^2_2}\omega^2
$$
and 
$$
\tau^2=\frac{1}{a^2_1}(F_0-F^2)\omega^1.
$$

\subsubsection{Curvatures}

We compute now $d\phi=\omega\wedge (A\omega^1+B\omega^2)$ and  
$
dw=C\omega\wedge\omega^1+ D\omega\wedge\omega^2+S\omega^1\wedge\omega^2 $(see Proposition \ref{Cartan-enriched})
.
From $d\phi=0$ we obtain $A=B=0$.  Computing $dw$ and comparing to the formula above we obtain
$$
C=\frac{F_1-F_{20}+FF_2}{3a^2_1a_2}, \ \ D=\frac{2F_2}{3a_1a^2_2}, \ \ S=\frac{F_{22}-F}{3a_1a_2}.
$$
In order to compute the global invariant we need to compute the coefficients $S_0,S_1$ and $S_2$ in equation \ref{dS} (a Bianchi identity) : $dS+2S\phi=S_0\omega+S_1\omega^1+S_2\omega^2$.  One obtains
$$
S_0=\frac{F_{220}-F_0}{3a^2_1a^2_2}, \ \ S_1=\frac{F_{221}-F_1}{3a^2_1a_2}, \ \ S_2=\frac{F_{222}-F_2}{3a_1a^2_2}.
$$

Now we use the expressions obtained in section \ref{section:embedddingY2-curvature} of the curvatures of $Y$ in terms of the curvature of $Y_2$.  In order to compute $\mu(Y)$ we need to compute $n$ and its derivatives (see
\ref{S0S1S2} and \ref{dn}).  

We have $2n=S_2+4D=\frac{F_{222}+7F_2}{3a^2_2a_1}$ and compute the left hand of $dn+3(\phi-w)n=n_0\omega+n_1\omega^1+n_2\omega^2$ (formula \ref{dn}) to obtain then
$$
n_0=\frac{1}{6a^2_1a^3_2}\left(F(F_{222}+F_2)+F_{2220}+7F_{20}\right),
$$
$$
n_1=\frac{1}{6a^2_1a^2_2}\left(F_2(F_{222}+F_2)+F_{2221}+7F_{21}\right)
$$
$$
n_2=\frac{1}{6a_1a^3_2}\left(F_{2222}+7F_{22}\right)
$$

We use equations \ref{n2} and \ref{p1} to compute the curvature functions $Q^1$ and $Q^2$.
We have $Q^1=-n_2+\frac 3 2 S\tau^1_2+\tau^1_{20}$ and $Q^2=-P_1-\frac 3 2 S\tau^2_1+\tau^2_{10}$.
For that sake, we compute first the derivatives of the torsion, $\tau^1_{20}$ and $\tau^2_{10}$, using formulas
\ref{dtau12} and \ref{dtau21}.
Computing the left hand side of the equation
$
d\tau^1_2+2\tau^1_2(\phi-3w)+(B+3D)\omega^1=\tau^1_{20}\omega+\tau^1_{22}\omega^2
$ and comparing to the right hand side we obtain
$$
\tau^1_{20}=-\frac{2F}{a_1a^3_2}.
$$

Analogously, computing the left hand side of the equation
$d\tau^2_1+2\tau^2_1(\phi+3w)-(A-3C)\omega^2=\tau^2_{10}\omega+\tau^2_{11}\omega^1$ we obtain
$$
\tau^2_{10}=\frac{F_{00}-4FF_0+2F^3}{a^3_1a_2}.
$$

\begin{prop}\label{proposition:Q1Q2DE}
Given a (local) differential equation as a path structure induced by the forms $\theta, \theta^1$ and $\theta^2$
as above one computes the curvature functions in terms of the enriched structure:
$$
Q^1=-\frac{1}{6a_1a^3_2}(F_{2222}+10F_{22}+9F)
$$
and
$$
Q^2=\frac{1}{6a^3_1a_2}(F_{2}F_{221}+3F_{2}F^2_1-4F_2F_{20}+4FF^2_2+F_{2211} +6F_{1}F_{11}
$$
$$
-4F_{201}+4F_{2}F_1+4FF_{21}-3F_{22}F_0+3F_{22}F^2-21FF_0+9F^3+6F_{00}).
$$
\end{prop}

\Pf Recall from formula \ref{S0S1S2} that $P= -S_1/2-2C$ and  we write $dP+3P(\phi+w)=P_0\omega+P_1\omega^1+P_2\omega^2$.  Computing the left side and comparing the right side
we  obtain the expression of $P_1$ which we use in the formulas above.
\EPf

The following proposition describes locally differential equations satisfying $Q^1=0$.

\begin{prop}\label{proposition:Q1=0}
Differential equations on an open subset with coordinates $(x,y)$ given by $
\theta=\cos{\alpha}dy-\sin{\alpha}dx,
$
$
\theta^1= \sin{\alpha}dy+\cos{\alpha}dx
$ 
and
$
\theta^2=d\alpha -F(x,y,\alpha)\theta^1
$ satisfy $$Q^1=0$$ if and only if
$$
F(x,y,\alpha)= A(x,y)\cos{\alpha}+B(x,y)\sin{\alpha}+C(x,y)\cos{3\alpha}+D(x,y)\sin{3\alpha}
$$
where $A,B,C$ and $D$ are functions on $x$ and $y$.
\end{prop}
\Pf
Observe that $Q^1=0$ is equivalent to $F_{2222}+10F_{22}+9F=F_{\alpha\alpha\alpha\alpha}+10F_{\alpha\alpha}+9F=0$. The only solutions to this linear equation are of the form above. 
\EPf

Using the coordinates $(x,y,p)$ as above where the differential equation is described as
$dp-G(x,y,p)dx=0$ the condition $Q^1=0$ implies that 
$G(x,y,p)$ is at most a third order polynomial in $p$ with coefficients functions of $x$ and $y$ (see \cite{A}).


\subsubsection{The global invariant}

We are ready now to use Proposition \ref{proposition-invariantY_2} to detemine the global invariant:
\begin{prop}\label{proposition;global-invariant}
Let $M$ be an enriched path structure defined by an ordinary differential equation of second order on the torus
with strict structure defined by the forms $\theta, \theta^1$ and $\theta^2$ as above.  Let $Y_2\rightarrow Y$ be the canonical embedding of the enriched geometry into the induced path geometry whose connection is $\pi$.  Then 
 $$
8\pi^2s^*(TC_2(\pi))=
\frac 1 {12}(-12F^2_\alpha +2(F_{\alpha\alpha\alpha x}\cos{\alpha}+F_{\alpha\alpha\alpha y}\sin{\alpha}+F_{\alpha\alpha\alpha\alpha}F)+14(F_{\alpha x}\cos{\alpha}+F_{\alpha y}\sin{\alpha}
$$
$$
+F_{\alpha\alpha}F)-3(-F_{\alpha \alpha x}\sin{\alpha}+F_{\alpha\alpha y}\cos{\alpha})
+3F_{\alpha}
-24(-F_x\sin\alpha+F_y\cos\alpha)+18F^2
+6F F_{\alpha\alpha})\theta\wedge \theta^1\wedge\theta^2.
$$
and
$$
8\pi^2\mu(Y)=8\pi^2\int_Ms^*(TC_2(\pi))=
\frac 1 {12}\int_M(-32F^2_\alpha +2F^2_{\alpha\alpha}+18F^2)\theta\wedge \theta^1\wedge\theta^2.
$$

 \end{prop}
 \Pf The terms in Proposition \ref{proposition-invariantY_2} were all computed before.  A substitution of these terms in the formula gives the first formula.  The second formula is obtained by integration by parts.
 
 \EPf

 \begin{cor}\label{corollary:Q1mu=0}
 Let $M$ be equipped with a path structure defined by an ordinary differential equation of second order on the torus
with strict structure defined by the forms $\theta, \theta^1$ and $\theta^2$.   Let $Y$ be the canonical Cartan bundle with its associated Cartan connection.
\begin{enumerate}
\item If $\theta^2=d\alpha-F(x,y)\theta^1$ (the function $F$ does not depend on $\alpha$).   Then $\mu(Y)=0$ if and only if $F=0$.
\item  $Q^1(Y)=0$ and $\mu(Y)=0$ if and only if $F=0$.
\end{enumerate}
 \end{cor}
 
 \Pf
 Clearly, if $F=0$ then $\mu(Y)=0$ and $Q^1=0$.  If $F$ does not depend on $\alpha$ then the invariant becomes 
 $$
 \frac{4}{3}\int_M F^2\theta\wedge \theta^1\wedge\theta^2,
$$
which is zero only if $F=0$.
Suppose now that  $Q^1(Y)=0$ and $\mu(Y)=0$.  Observe that the integral formula for the invariant, by an integration by part and a slight rearrangement, may be written as
$$
\frac 1 {12}\int_M(-12F^2_\alpha +2F_{\alpha\alpha\alpha\alpha}F+20F_{\alpha\alpha}F +18F^2)\theta\wedge \theta^1\wedge\theta^2.
$$
Using the expression of $Q^1$ given in \ref{proposition:Q1Q2DE} and the hypothesis $Q^1=0$, we obtain 
$
F_{2222}+10F_{22}+9F=0$, and therefore
$$
2F_{\alpha\alpha\alpha\alpha}F+20F_{\alpha\alpha}F +18F^2=0.
$$
Therefore
$$
8\pi^2\mu(Y)=-\int_MF^2_\alpha \theta\wedge \theta^1\wedge\theta^2.
$$ 
We observe therefore that if $\mu(Y)=0$ then $F_\alpha$ should be null.  But if $F$ does not depend on $\alpha$ it should be null by the first part.
 \EPf
 

\section{Path structures on a torus}
We recall example 
{\bf{II}} which is the torus $T^3$ with coordinates $(x,y,t)$ ($\mod 1$) and the global contact form, for a fixed $n\in \Z^*$,
$$
\theta= \cos (2\pi n t) dx- \sin(2\pi n t) dy.
$$

It was proven independently by E. Giroux and Y. Kanda that the contact structures defined by these contact forms classify all tight structures
on $T^3$ (see \cite{Y}).  We will show here that for each of these contact structures one can define  a flat path structure.

There are two canonical global vector fields on the distribution given by
$X_1=\frac{\partial}{\partial t}$ and $X_2=\sin (2\pi n t) \frac{\partial}{\partial x}+\cos(2\pi n t)\frac{\partial}{\partial y}$.

We define
$$
 \theta^1=-2\pi n dt,\ \  \theta^2= \sin (2\pi n t) dx+\cos(2\pi n t) dy,
$$
so that $d\theta =\theta^1\wedge \theta^2$ and we define the strict path structure defined by these forms.
We compute
$$
d\theta^1=0, \ \ d\theta^2=-\theta^1\wedge \theta.
$$
Comparing now with the enriched path connection we obtain
$$
\phi=\frac 1 2(\frac{da_1}{a_1}+\frac{da_2}{a_2}),\ \ 3w=\frac 1 2(\frac{da_1}{a_1}-\frac{da_2}{a_2}),\ \ \tau^1_2=0, \ \ \tau^2_1=-\frac{1}{a_1^2}
$$
and therefore $d\phi=d w=d\tau^1_2=d\tau_1^2+2\tau^2_1(\phi+3w)=0$. 
It follows that $A=B=C=D=S=\tau^1_{20}=\tau^2_{10}=0$, and it follows from formulas \ref{n2} and \ref{p1} that $Q^1=Q^2=0$.
We proved:
\begin{lemma}  
The path structures defined by the forms $\theta^1,\theta^2,\theta$ on $T^3$ are flat.
\end{lemma}

We now define a new strict path structure by fixing the contact form $\theta$ and changing $\theta^1$ and $\theta^2$ by a constant matrix:
$$
\left ( \begin{array}{c}

                       {\theta^1}' \\

                      { \theta^2}'                      

                \end{array} \right )
=\left ( \begin{array}{cc}

                        a      &    b  \\

                        c  	&      f        \\

                \end{array} \right ) \left ( \begin{array}{c}

                       \theta^1 \\

                       \theta^2                      

                \end{array} \right )  \ \  \ \ \mbox{or}\ \ 
              \left ( \begin{array}{c}

                       {\theta^1}\\

                      { \theta^2}                      

                \end{array} \right )
=\left ( \begin{array}{cc}

                        f      &   - b  \\

                       - c  	&      a      \\

                \end{array} \right ) \left ( \begin{array}{c}

                       {\theta^1}' \\

                       {\theta^2}'                      

                \end{array} \right ),
$$
where 
$$
\det \left ( \begin{array}{cc}

                        a      &    b  \\

                        c  	&      f        \\

                \end{array} \right )=1.
$$
We compute $d{\theta^1}'= ad{\theta^1}+bd{\theta^2}=-b\theta^1\wedge \theta=-bf{\theta^1}'\wedge \theta+b^2{\theta^2}'\wedge \theta$
and $d{\theta^2}'= cd{\theta^1}+fd{\theta^2}=-f\theta^1\wedge \theta=-f^2{\theta^1}' \wedge \theta+bf{\theta^2}'\wedge   \theta$.
In order to compute the  enriched connection we need to find $\phi', w', {\tau^1}', {\tau^2}'$ satisfying
$$
d{\omega^1}'=(\phi'+3w')\wedge{\omega^1}'+\omega' \wedge {\tau^1_2}'{\omega^2}', \ \ \ d{\omega^2}'=(\phi'-3w')\wedge{\omega^2}'-\omega' \wedge {\tau^2_1}'{\omega^1}'.
$$Comparing with the structure equations and observing ${\omega^i}'=a_i{\theta^i}'$ we obtain
$$
\phi'=\frac 1 2(\frac{da_1}{a_1}+\frac{da_2}{a_2}),\ \ 3w'=\frac 1 2(\frac{da_1}{a_1}-\frac{da_2}{a_2})+\frac{bf}{a_1a_2}\omega,\ \ 
$$
$$
\ {\tau^1_2}'=-\frac{b^2}{a_2^2}{\theta^2}'\ \ \mbox{and}\ \   {\tau^2_1}'=-\frac{f^2}{a_1^2}{\theta^1}'.$$
Then $d\phi=0, \ \ 3dw=\frac{bf}{a_1a_2}{\omega^1}'\wedge{\omega^2}'$, and it follows $A=B=C=D=0$ and $3S=\frac{bf}{a_1a_2}$. Also
$$
d{\tau^1_{2}}'=-2{\tau^1_{2}}'(\phi-3w)+\frac{2b^3f}{a_2^3a_1}\omega,\ \ d{\tau^2_{1}}'=-2{\tau^2_{1}}'(\phi+3w)-\frac{2bf^3}{a_1^3a_2}\omega,
$$
and ${\tau^1_{20}}'=\frac{2b^3f}{a_2^3a_1}, \ \ {\tau^2_{10}}'=-\frac{2bf^3}{a_1^3a_2}$.
At last $dS=-2S\phi$ and we get $S_0=S_1=S_2=0$. Then $P=n=0$, and we obtain from formulas \ref{n2} and \ref{p1} that 
$$
Q^1=\frac{3}{2}\frac{b^3f}{a_1a_2^3}, \ \ Q^2=-\frac{3}{2}\frac{bf^3}{a_1^3a_2}.
$$

We proved
\begin{lemma}  
The path structures defined by the forms $\theta'^1,\theta'^2,\theta$ on $T^3$ have curvatures $Q^1=\frac 3 2 b^3f, \ \ \ Q^2=-\frac 3 2 bf^3$ (computed through a section on the torus).
\end{lemma}
 
 Note that the path structure is flat if and only if the one of the torsions ${\tau^1}'$ or ${\tau^2}'$ are zero and this happens if the direction defined by $\frac{\partial}{\partial t}$ is one of the line bundles contained in the contact bundle of the path structure. The couple $(b,f)$ is determined up to a sign by the curvatures
 $Q^1$ and $Q^2$.
 
 The global invariant is given in the next Proposition.
 
 \begin{prop}\label{proposition:torustight}
Let ${T^3}_n(a,b,c,d,f)$ as the path structure on the torus defined as above.  Then the global invariant is
$$
\mu({T^3}_n(a,b,c,d,f))=\frac{3n}{8\pi}(bf)^2.
$$
 \end{prop}
 
 \Pf  This is a direct computation using the formula for the global invariant (see formula \ref{pi3}):
 $$
\int_{T^3}s^*TC_2(\pi)=\int_{T^3}\frac{1}{8\pi^2} (2{\tau^1_{2}}'{\tau^2_{1}}'\theta\wedge \theta^1\wedge \theta^2-\frac 9 2 w\wedge \theta^1\wedge \theta^2)=\int_{T^3}\frac{1}{8\pi^2} \frac{3}{2}b^2f^2\theta\wedge \theta^1\wedge \theta^2.
$$
Therefore
 $$
\mu({T^3}_n(a,b,c,d,f))=\int_{T^3} \frac{3}{16\pi^2}(bf)^2\theta d\theta=\frac{3n}{8\pi}(bf)^2.
$$
 
 \EPf
 
 Note that the global invariant is null if and only if the path structure is flat.

\section{Invariant path structures on ${ \SU}(2)$}

Tight contact structures on  $S^3$ are all contactomorphic (see \cite{E}).    In this section we explicit homogeneous strict path structures on ${ \SU}(2)$ which are carried by a fixed left invariant tight contact structure.

Let $\alpha, \beta, \gamma$ be a basis of left invariant 1-forms defined
on $\SU(2)$ with
$$
d\alpha= -\beta\wedge \gamma,\ \ \ d\beta= -\gamma\wedge \alpha,\ \ \ d\gamma= -\alpha\wedge \beta\ \ \ 
$$
A strict path structure on ${ \SU}(2)$ is given by fixing the contact form $\theta= \gamma$ and 
the line fields $E^1=\ker{\alpha}\cap \ker{\theta}$ and $E^2=\ker{\beta}\cap \ker{\theta}$.

We define strict path structures by choosing a map from $\SU(2)$ to $\SL(2,\R)$:
$$
\theta=\gamma, \ \ \ Z^1=r_1\beta+r_2\alpha, \ \ \ Z^2=s_1\beta+s_2\alpha,
$$
with $r_1s_2-r_2s_1=1$.
Then
$$
d\theta=Z^1\wedge Z^2.
$$
 
In the case the map $\SU(2)\rightarrow \SL(2,\R)$ is constant, from 
$\beta =s_2Z^1-r_2Z^2$ and $\alpha =-s_1Z^1+r_1Z^2$,
we obtain
$$
dZ^1=r_1d\beta +r_2 d\alpha=\theta \wedge \left( xZ^1+yZ^2\right)
$$
and analogously,
$$
dZ^2=\theta \wedge \left( zZ^1-xZ^2\right),
$$
 where
$$
x=r_1s_1+r_2s_2,\ \ \ y=-(r_1^2+r_2^2), \ \ \ z=s_1^2+s_2^2.
$$
Observe that $x^2+yz=-1$.
Then for a enriched path structure with coframes obtained from t!he tautological forms $\omega=a_1a_2\theta$, $\omega^1=a_1Z^1$ and $\omega^2=a_2Z^2$ we obtain
$$
d\omega^1=(\frac{da_1}{a_1}+x\theta)\wedge \omega^1+a_1y\theta\wedge Z^2.
$$
$$
d\omega^2=(\frac{da_2}{a_2}-x\theta)\wedge \omega^2+a_2z\theta\wedge Z^1
$$
From Proposition \ref{cartan-strict} we have
$$
\phi=\frac 1 2(\frac{da_1}{a_1}+\frac{da_2}{a_2}),\ \ 3w=\frac 1 2(\frac{da_1}{a_1}-\frac{da_2}{a_2})+\frac{x}{a_1a_2}\omega,\ \ 
$$
$$
\tau^1_2=\frac{y}{ a_2^2} , \ \ \ \tau^2_1=-\frac{z}{ a_1^2}.
$$
and therefore
$$
d\phi=0, \ \  3dw=d(x\theta)= \frac{x}{a_1a_2}\omega^1\wedge \omega^2
$$
so that $S=\frac{x}{3a_1a_2}, A=B=C=D=0$.

From 
$$
d\tau_2^1=-2\frac{da_2}{a_2}\frac{y}{a_2^2}=-2\tau^1_2(\phi-3w)-2\frac{xy}{a_1a_2^3}\omega
$$
$$
d\tau_1^2=2\frac{da_1}{a_1}\frac{z}{a_1^2}=-2\tau^2_1(\phi+3w)-2\frac{xz}{a_2a_1^3}\omega
$$
we obtain 
$$\tau_{20}^1=-2\frac{xy}{a_1a_2^3}, \ \  \tau_{10}^2=-2\frac{xz}{a_2a_1^3}$$
From
$$
dS=-\frac{x}{3a_1a_2}(\frac{da_1}{a_1}+\frac{da_2}{a_2})=-2\phi S
$$
we obtain $S_0=S_1=S_2=0$, ans $P=n=0$.

It follows from formulas \ref{n2} and \ref{p1} that $Q^1=\tau^1_{20}+\frac 3 2 S\tau^1_2$ and $Q^2=\tau^2_{10}-\frac 3 2 S\tau^2_1, $
therefore
  $$
Q^1=-\frac{xy}{a_1a_2^3}
$$
and
$$
Q^2=-\frac{xz}{a_1^3a_2}.
$$
Observe that $y$ and $z$ never vanish.  We conclude that the invariant strict structure on $\SU(2)$ is a flat path structure 
if and only if $x=0$.  This  can be interpreted, because $x=r_1s_1+r_2s_2$, as the strict structures such that the directions $E^1$ and $E^2$ are perpendicular for the canonical metric defined by the forms $\alpha$ and $\beta$.

\begin{prop}\label{proposition:strictori}
Define strict path structures on $\SU(2)$ by choosing a constant map from $\SU(2)$ to $\SL(2,\R)$:
$$
\theta=\gamma, \ \ \ Z^1=r_1\beta+r_2\alpha, \ \ \ Z^2=s_1\beta+s_2\alpha,
$$
with $r_1s_2-r_2s_1=1$.  Let $x=r_1s_1+r_2s_2$. Then the global invariant of the induced path structure is
$$
\mu(\SU(2)(r_1,r_2,s_1,s_2))=-\frac{1}{2}-\frac{3}{8}x^2.
$$
\end{prop}
 \Pf
We compute, using formula \ref{pi3}, the global invariant for the family of structures defined on $\SU(2)$.
We have from above that $
x=r_1s_1+r_2s_2, y=-(r_1^2+r_2^2), z=s_1^2+s_2^2
$ and that $x^2+yz=-1$.  Then it follows
$$
\int_{\SU(2)}s^*TC_2(\pi)=\int_{\SU(2)}\frac{1}{8\pi^2} (2{\tau^1_{2}}{\tau^2_{1}}\theta-\frac 9 2 S w)\wedge \theta^1\wedge \theta^2=-\int_{\SU(2)}\frac{1}{8\pi^2}(2yz+\frac{1}{2}x^2)\gamma\wedge \beta\wedge \alpha
$$
$$
=
\int_{\SU(2)} \frac{1}{8\pi^2}(-2-\frac{3}{2}x^2)\gamma\wedge \beta\wedge \alpha.
$$
We use then that $\int_{\SU(2)} \gamma\wedge \beta\wedge \alpha=2\pi^2$.
\EPf

Observe that the invariant is never null for this family even in the case of a flat path structure (which happens when $x=0$). Also the critical point of the invariant along this family is a maximal at $x=0$, at a flat structure, and it is equal to $-\frac{1}{2}$.


\begin{flushleft}
  \textsc{E. Falbel\\
  Institut de Math\'ematiques \\
  de Jussieu-Paris Rive Gauche \\
CNRS UMR 7586 and INRIA EPI-OURAGAN \\
 Sorbonne Universit\'e, Facult\'e des Sciences \\
4, place Jussieu 75252 Paris Cedex 05, France \\}
 \verb|elisha.falbel@imj-prg.fr|
 \end{flushleft}
\begin{flushleft}
  \textsc{J. M.  Veloso\\
  Faculdade de Matem\' atica - ICEN\\
Universidade Federal do Par\'a\\66059 - Bel\' em- PA - Brazil}\\
  \verb|veloso@ufpa.br|
\end{flushleft}

\end{document}